\documentclass[aop,seceqn,MSNbibl,citesort,dvips]{arximspdf}

\doi{10.1214/09-AOP523}
\volume{38}
\issue{5}
\pubyear{2010}
\firstpage{1817}
\lastpage{1869}

\makeatletter

\newtheorem{theorem}{Theorem}[section]
\newtheorem{prop}[theorem]{Proposition}
\newtheorem{corollary}[theorem]{Corollary}
\newtheorem{lemma}[theorem]{Lemma}
\newproclaim{remark}[theorem]{Remark}
\newproclaim{example}[theorem]{Example}

\makeatother

\begin{document}
\begin{frontmatter}

\title{A change of variable formula with It\^o~correction~term}
\runtitle{Change of variable with It\^o term}

\begin{aug}
\author[A]{\fnms{Krzysztof} \snm{Burdzy}\thanksref{t1}\ead[label=e1]{burdzy@math.washington.edu}} and
\author[B]{\fnms{Jason} \snm{Swanson}\corref{}\thanksref{t2}\ead[label=e2]{swanson@mail.ucf.edu}}
\runauthor{K. Burdzy and J. Swanson}
\affiliation{University of Washington and University of Central Florida}
\address[A]{Department of Mathematics\\
University of Washington\\
Box 354350\\
Seattle, Washington 98195-4350\\
USA\\
\printead{e1}} 
\address[B]{Department of Mathematics\\
University of Central Florida\\
4000 Central Florida Blvd.\\
P.O. Box 161364\\
Orlando, Florida 32816-1364\\
USA\\
\printead{e2}}
\end{aug}

\thankstext{t1}{Supported in part by NSF Grant DMS-09-06743 and by
Grant N N201 397137, MNiSW, Poland.}

\thankstext{t2}{Supported in part by the VIGRE grant of the University
of Wisconsin-Madison and by NSA Grant H98230-09-1-0079.}

\received{\smonth{4} \syear{2008}}
\revised{\smonth{8} \syear{2009}}

%
\begin{abstract}
We consider the solution $u(x,t)$ to a stochastic heat equation.
For fixed~$x$, the process $F(t)=u(x,t)$ has a nontrivial quartic
variation. It follows that $F$ is not a semimartingale, so a
stochastic integral with respect to $F$ cannot be defined in the
classical It\^o sense. We show that for sufficiently
differentiable functions $g(x,t)$, a stochastic integral $\int
g(F(t),t) \,dF(t)$ exists as a limit of discrete, midpoint-style
Riemann sums, where the limit is taken in distribution in the
Skorokhod space of cadlag functions. Moreover, we show that this
integral satisfies a change of variable formula with a
correction term that is an ordinary It\^o integral with respect to
a Brownian motion that is independent of $F$.
\end{abstract}

%

\end{frontmatter}

\section{Introduction}\label{S:intro}

Recall that the classical It\^o formula (i.e., change of variable
formula) contains a ``stochastic correction term'' that is a
Riemann integral. A purely intuitive conjecture is that the It\^o
integral itself may appear as a stochastic correction term in a
change of variable formula when the underlying stochastic process
has fourth order scaling properties. The first formula of this
type was proven in \cite{BM}; however, the ``fourth order
scaling'' process considered in that paper was a highly abstract
object with little intuitive appeal. The present article presents
a change of variable formula with It\^o correction term for a
family of processes with fourth order local scaling properties;
see (\ref{newito}) and Corollary \ref{C:main}.

The process which is our primary focus is the solution, $u(x,t)$,
to the stochastic heat equation $\partial_t u=\frac12 \,\partial_x^2
u+\dot{W}(x,t)$ with initial conditions $u(x,0)\equiv0$, where
$\dot{W}$ is a space--time white noise on $\mathbb{R}\times[0,\infty)$.
That is,
%
\begin{equation}\label{SPDE}
u(x,t) = \int_{\mathbb{R}\times[0,t]} p(x - y,t - r)W(dy\times dr),
\end{equation}
where $p(x,t)=(2\pi t)^{-1/2}e^{-x^2/2t}$ is the heat kernel. Let
$F(t)=u(x,t)$, where $x\in\mathbb{R}$ is fixed. In the prequel to this
paper \cite{S}, it was shown that $F$ is a continuous, centered
Gaussian process with covariance function
%
\begin{equation}
\rho(s,t) = EF(s)F(t) = (2\pi)^{-1/2}(|t + s|^{1/2} - |t - s|^{1/2})
\end{equation}
and that $F$ has a nontrivial quartic variation. In particular,
\[
\sum_{j=1}^n \bigl|F(j/n)-F\bigl((j-1)/n\bigr)\bigr|^4 \to\frac6\pi
\]
in $L^2$. It follows that $F$ is not a semimartingale, so a
stochastic integral with respect to $F$ cannot be defined in the
classical It\^o sense. In this paper, we complete the construction
of a stochastic integral with respect to $F$ which is a limit of
discrete Riemann sums.

More generally, we shall construct a stochastic integral with respect
to any process $X$ of the form $X=cF+\xi$, where $c\in\mathbb{R}$
and $\xi
$ is a
continuous
stochastic process, independent of $F$, satisfying
%
\begin{equation}\label{exten0}
\xi\in C^1((0,\infty)) \quad\mbox{and}\quad
\overline{\lim_{t\to0}} |\xi'(t)| < \infty\qquad\mbox{a.s.}
\end{equation}
This allows us, for example, to consider solutions to (\ref{SPDE})
with nonzero initial conditions. Another example of such an $X$
is fractional Brownian motion with Hurst parameter $1/4$; see
Examples \ref{ex:1} and \ref{ex:2} for more details.

Note that $\xi$ (and therefore $X$) need not be a Gaussian process. If
it is Gaussian, however, its mean function will be $\mu_X(t) = EX(t) =
\mu_\xi(t)$ and its covariance function will be $\rho_X(s,t)
=c^2\rho
(s,t) + \rho_\xi(s,t)$. Conversely, the results in this paper will
apply to any Gaussian process $X$ whose mean and covariance have the
respective forms $\mu_X=\widetilde\mu$ and $\rho_X =c^2\rho
+\widetilde\rho$, where
$\widetilde\mu$ and $\widetilde\rho$ are the mean and covariance,
respectively,
of a
Gaussian process satisfying (\ref{exten0}).

We conjecture that the results in this paper hold when $\xi$ is only
required to be of bounded variation. We require $\xi$ to be $C^1$,
however, because of our particular method of proof; see the proofs of
Corollaries 4.6 and 6.4 for further details.

For simplicity, we consider only evenly spaced partitions. That is,
given a positive integer $n$, let $\Delta t=n^{-1}$, $t_j=j\Delta t$
and $\Delta
X_j=X (t_j)-X(t_{j-1})$. Let $\lfloor x\rfloor$ denote the greatest
integer less
than or equal to $x$. For $g\in C(\mathbb{R}\times[0,\infty))$, we consider
the midpoint-style Riemann sums
%
\begin{equation}\label{riem}
I_n^X(g,t) = \sum_{j=1}^{\lfloor nt/2\rfloor}g(X(t_{2j-1}),t_{2j-1})
\bigl(X(t_{2j}) - X(t_{2j-2})\bigr).
\end{equation}
When $X=F$, we will simply write $I_n$, rather than $I_n^F$.

In the construction of the classical It\^o integral, the quadratic
variation of the integrator plays a crucial role. Although the
quadratic variation of $X$ is infinite, the ``alternating quadratic
variation'' of $X$ is finite. That is, $Q_n^X(t)=\sum_{j=1}^{\lfloor
nt/2\rfloor}(\Delta X_{2j}^2-\Delta X_{2j-1}^2)$ converges in law.
If we denote
the limit process by $\{X\}_t$, then it is a simple corollary of the
main result in \cite{S} that $\{X\}_t$ is a Brownian motion which is
independent of $X$. More specifically, $(X,Q_n^X)\to(X,\kappa c^2B)$,
where $B$ is a standard Brownian motion, independent of $X$, and
$\kappa
\approx1.029$ [see (\ref{kap}) for the precise definition of $\kappa$].
The convergence here is in law in the Skorokhod space of cadlag
functions from $[0,\infty)$ to $\mathbb{R}^2$, denoted by $D_{\mathbb{R}
^2}[0,\infty)$.

We shall show that $I_n^X(g,t)$ also converges in law. If $\int_
0^t g(X(s),s) \,dX(s)$ denotes a process with this limiting law,
then our main result (Corollary \ref{C:main}) is the following
change of variable formula:
\begin{eqnarray*}
g(X(t),t) &=& g(X(0),0) + \int_0^t \partial_t g(X(s),s) \,ds
+ \int_0^t \partial_x g(X(s),s) \,dX(s)\\
&&{} + \frac12 \int_0^t \partial_x^2 g(X(s),s) \,d\{X\}_s,
\end{eqnarray*}
where the equality is in law as processes. This can be rewritten as
%
\begin{eqnarray}\label{newito}\qquad
g(X(t),t) &=& g(X(0),0) + \int_0^t \partial_t g(X(s),s) \,ds
+ \int_0^t \partial_x g(X(s),s) \,dX(s)\nonumber\\[-8pt]\\[-8pt]
&&{} + \frac{\kappa c^2}2 \int_0^t \partial_x^2 g(X(s),s) \,dB(s),\nonumber
\end{eqnarray}
where this last integral is a classical It\^o integral with respect to
a standard Brownian motion that is independent of $X$.

To state our results more completely, let $Y$ be a semimartingale and define
%
\begin{eqnarray}\label{formint}
I^{X,Y}(\partial_x g,t) &=& g(X(t),t) - g(X(0),0)
- \int_0^t \partial_t g(X(s),s) \,ds\nonumber\\[-8pt]\\[-8pt]
&&{} - \frac{\kappa}2\int_0^t \partial_x^2 g(X(s),s) \,dY(s).\nonumber
\end{eqnarray}
We show that
\[
(F, Q_n^F, I_n^X(\partial_x g,\cdot)) \to(F, \kappa B,
I^{X,c^2B}(\partial_x
g,\cdot))
\]
in law in $D_{\mathbb{R}^3}[0,\infty)$
whenever $g\in C^{9,1}_4(\mathbb{R}\times[0,\infty))$. [See
(\ref{c1k0})--(\ref{c1k3}) for the precise definition of the space
$C^{k,1}_r$. Also see Remarks \ref{R:main1} and \ref{R:main2}.]

The benefit of having the convergence of this triple, rather than
just the Riemann sums, can be seen if one considers two separate
sequences of sums: $\{I_n^{X_1}(g_1, \cdot)\}$ and $\{I_n^{X_2}
(g_2,\cdot)\}$. As $n\to\infty$, these sequences will converge
jointly in law. Separately, each limit will satisfy
(\ref{newito}); moreover, the Brownian motions which appear in
the two limits will be identical. In this sense, the Brownian
motion in (\ref{newito}) depends only on $F$ and not on $\xi$,
$c$ or $g$. Clearly, this can be extended to any finite
collection of sequences of Riemann sums.

In the course of our analysis, we will also obtain the asymptotic
behavior of the trapezoid-style sum
%
\begin{equation}\label{Triem}
T_n^X(g,t) = \sum_{j=1}^{\lfloor nt\rfloor}\frac{
g(X(t_{j-1}),t_{j-1}) + g(X(t_j),t_j)}2\Delta X_j.
\end{equation}
We shall see (Corollary \ref{C:expan3}) that $T_n^X(\partial_x
g,t)\to
g(X(t),t)-g(X(0),0) - \int_0^t \partial_tg(X(s),s) \,ds$ uniformly on
compacts in probability (ucp) whenever $g\in
C^{7,1}_3(\mathbb{R}\times[0,\infty))$. This result remains true
even when
$X=cF+\xi$, where $\xi$ satisfies only (\ref{exten0}), and is not
necessarily independent of $F$.

It is instructive to contrast these results with those of Russo,
Vallois and coauthors \cite{G1,G2,R1,R2}, who, in the context of
fractional Brownian motion, use a regularization procedure to transform
these Riemann sums into integrals before passing to the limit; see also
\cite{C}. For instance, if $g$ does not depend on $t$, then the
regularized midpoint sum is
\[
\frac1{2\varepsilon}\int_0^t g'(F(s))\bigl(F(s+\varepsilon) -
F\bigl((s-\varepsilon)\vee0\bigr)\bigr) \,ds
\]
and the regularized trapezoid sum is
\[
\frac1{2\varepsilon}\int_0^t \bigl(g'(F(s)) + g'\bigl(F(s+\varepsilon
)\bigr)\bigr)\bigl(F(s+\varepsilon) - F(s)\bigr) \,ds.
\]
Using a change of variables, we can see that if $g'$ is locally
integrable, then the difference between these two integrals goes to
zero almost surely as $\varepsilon\to0$. Hence, under the regularization
procedure, the
midpoint and trapezoid sums exhibit the same limiting behavior: they
converge ucp to integrals satisfying the classical change of variable
formula from ordinary calculus. Under the discrete approach which we
are following, however, we see new behavior for the midpoint sum: the
emergence of a
correction term which is a classical It\^o integral against an
independent Brownian motion.

It should be noted that all of our convergence results rely on the fact
that $F$ is a quartic variation process.
That is,
%
\begin{equation}\label{order}
C_1\Delta t^{2H} \le E\Delta F_j^2 \le C_2\Delta t^{2H},
\end{equation}
where $H=1/4$. For example, the convergence of $Q_n^F$ to a
Brownian motion is made plausible by the fact that it is a sum of
terms of the form $\Delta F_{2j}^2 - \Delta F_{2j-1}^2$, each of which
is approximately mean zero with an approximate variance of $\Delta
t$. If we replace $F$ by a rougher process which satisfies
(\ref{order}) for some $H<1/4$, then the midpoint sums will
evidently diverge. On the other hand, the ucp convergence of the
trapezoid sums $T_n(\partial_x g,t)$ remains plausible
for any $H>1/6$. This is consistent with the analogous results in
\cite{C,G1} for regularized sums.

The critical case for the trapezoid sum is $H=1/6$. At the time of
writing, we know of only one result in this case. If
$g(x,t)=x^3$, then
\[
T_n(\partial_x g,t) \approx F(t)^3 - F(0)^3
+ \frac12\sum_{j=1}^{\lfloor nt\rfloor}\Delta F_j^3.
\]
[Here, and in what follows, $X_n(t)\approx Y_n(t)$ shall mean that $X_n
- Y_n \to0$ ucp.] If $F$ is replaced by fractional Brownian motion
with Hurst parameter $H=1/6$, then this last sum converges in law to a
Brownian motion; see \cite{N}, for example. It is natural to conjecture
that a result analogous to (\ref{newito}) also holds in this case.

Our project is related to, and inspired by, several areas of stochastic
analysis. Recently, a new approach to integration was developed by
T. Lyons (with coauthors, students and other researchers). The new method is
known as ``rough paths''; an introduction can be found in \cite{LCL}.
Our approach is much more elementary since it is based on a form of
Riemann sums. We consider it of interest to see how far the classical
methods can be pushed and what they can yield. The It\^o-type
correction term in our change of variable formula has a certain
elegance to it, and a certain logic, if we recall that our underlying
process has quartic variation. Finally, our project can be considered a
toy model for some numerical schemes. The fact that the correction term
in the change of variable formula involves an independent Brownian
motion may give some information about the form and size of errors in
numerical schemes.

After the first draft of this paper had been finished, we received
a preprint \cite{NR} from Nourdin and R\'eveillac,
prepared independently of ours and using different methods. That
paper contains a number of results, one of which, Theorem 1.2, is
a special case of our Corollary \ref{C:main}. Namely, if
$X=B^{1/4}$ (fractional Brownian motion with Hurst parameter $H
=1/4$), if $g$ does not depend on $t$ and if $g$ satisfies an
additional moment condition (see $\mathbf H_q$ in Section 3 of
\cite{NR}), then \cite{NR} gives the
convergence in distribution of the scalar-valued random variables
$I_n^X(g',1)$. While \cite{NR} is devoted exclusively to
fractional Brownian motion, it is mentioned in a footnote that a
Girsanov-type transformation can be used to extend the results
from $B^{1/4}$ to $F$.

\eject
\section{Preliminaries}

\subsection{Tools for cadlag processes}

Here,  and in the remainder of this paper, $C$~shall denote a
constant whose value may change from line to line.

Let $D_{\mathbb{R}^d}[0,\infty)$ denote the space of cadlag functions from
$[0, \infty)$ to $\mathbb{R}^d$ endowed with the Skorokhod topology.
We use
the notation $x(t-)=\lim_{s\uparrow t}x(s)$ and $\Delta x(t)=x(t)-x(t-)$.
Note that if $F_n(t) =F(\lfloor nt\rfloor/n)$, then $\Delta
F_n(t_j)=F(t_j)-F(t_{j-1})$. As in Section \ref{S:intro}, we shall
typically use $\Delta F_j$ as a shorthand notation for $\Delta F_n(t_j)$.

We note for future reference that if $x$ is continuous, then $x_n\to x$
in the Skorokhod topology if and only if $x_n\to x$ uniformly on
compacts. For our convergence results, we shall use the following
moment condition for relative compactness, which is a consequence of
Theorem 3.8.8 in \cite{EK}.
\begin{theorem}\label{T:momcrit}
Let $\{X_n\}$ be a sequence of processes in $D_{\mathbb{R}^d}[0,\infty)$.
Let $q(x) = |x|\wedge1$. Suppose that for each $T>0$, there exist
$\nu>0$, $\beta>0$, $C>0$ and $\theta>1$ such that:
\begin{longlist}
\item$E[q(X_n(t+h)-X_n(t))^{\beta/2}q(X_n(t)-X_n(t-h))^{\beta/2}]
\le Ch^\theta$ for all $n$ and all $0\le t\le T+1$, $0\le h\le t$;
\item$\lim_{\delta\to0}\sup_n E[q(X_n(\delta)-X_n(0))^\beta]=0$;
\item$\sup_n E[|X_n(T)|^\nu]<\infty$.
\end{longlist}
Then $\{X_n\}$ is relatively compact, that is, the distributions are
relatively compact in the topology of weak convergence.
\end{theorem}
\begin{corollary}\label{C:momcrit}
Let $\{X_n\}$ be a sequence of processes in $D_{\mathbb{R}^d}[0,\infty
)$. Let
$q(x) = |x|\wedge1$. Let $\varphi_1,\varphi_2$ be nonnegative
functions of $n$
such that $\sup_nn^{-1} \varphi_1(n)\times\break\varphi_2(n)<\infty$. Suppose
that for each $T>0$, there exist $\nu>0$, $\beta>0$, $C>0$ and $\theta>1$ such that
$\sup
_n E[|X_n(T)|^\nu]<\infty$ and
%
\begin{equation}\label{momcrit}
E\bigl[q\bigl(X_n(t) - X_n(s)\bigr)^\beta\bigr]
\le C \biggl(\frac
{\varphi_2(n)\lfloor\varphi_1(n)t\rfloor- \varphi_2(n)\lfloor
\varphi_1(n)s\rfloor}n \biggr)^\theta
\end{equation}
for all $n$ and all $0\le s,t\le T$. Then $\{X_n\}$ is relatively compact.
\end{corollary}
\begin{pf}
We apply Theorem \ref{T:momcrit}. By hypothesis, condition (iii)
holds. Taking $s=0$ and $t=\delta$ in (\ref{momcrit}) gives condition
(ii). By H\"older's inequality,
\begin{eqnarray*}
&&E\bigl[q\bigl(X_n(t+h)-X_n(t)\bigr)^{\beta/2}q\bigl(X_n(t)-X_n(t-h)\bigr)^{\beta/2}\bigr]\\
&&\qquad\le C \biggl(\frac
{\varphi_2(n)\lfloor\varphi_1(n)(t+h)\rfloor- \varphi
_2(n)\lfloor\varphi_1(n)t\rfloor}n
\biggr)^{\theta/2}\\
&&\qquad\quad{}\times\biggl(\frac
{\varphi_2(n)\lfloor\varphi_1(n)t\rfloor- \varphi_2(n)\lfloor
\varphi_1(n)(t-h)\rfloor}n
\biggr)^{\theta/2}.
\end{eqnarray*}
If $\varphi_1(n)h<1/2$, then the right-hand side of the above
inequality is
zero. Assume that $\varphi_1(n)h\ge1/2$. Then
\begin{eqnarray*}
&&
E\bigl[q\bigl(X_n(t+h)-X_n(t)\bigr)^{\beta/2}q\bigl(X_n(t)-X_n(t-h)\bigr)^{\beta/2}\bigr]\\
&&\qquad\le C \biggl(\frac{\varphi_2(n)\varphi_1(n)h + \varphi
_2(n)}n \biggr)^\theta
\le\widetilde C \biggl(h + \frac1{\varphi_1(n)} \biggr)^\theta
\le\widetilde C(3h)^\theta,
\end{eqnarray*}
which verifies condition (i).
\end{pf}

In general, the relative compactness in $D_\mathbb{R}[0,\infty)$ of
$\{X_n\}$ and $\{Y_n\}$ does not imply the relative compactness of
$\{X_n+Y_n\}$. This is because addition is not a continuous
operation from $D_\mathbb{R}[0,\infty)^2$ to $D_\mathbb{R}[0,\infty
)$. It is,
however, a continuous operation from $D_{\mathbb{R}^2} [0,\infty)$ to
$D_{\mathbb{R}}[0,\infty)$. To make use of this, we shall need the
following well-known result and its subsequent corollary.
\begin{lemma}\label{L:conlim}
Suppose that $x_n\to x$ in $D_{\mathbb{R}}[0,\infty)$ and $y_n\to y$ in
$D_{\mathbb{R}
}[0, \infty)$. If $\Delta x(t)\Delta y(t)=0$ for all $t\ge0$, then
$x_n+y_n\to x+y$ in $D_ {\mathbb{R}} [0,\infty)$.
\end{lemma}
\begin{corollary}\label{C:conlim2}
Suppose that the sequences $\{X_n\}$ and $\{Y_n\}$ are relatively
compact in $D_{\mathbb{R}}[0,\infty)$. If every subsequential limit
of $\{
Y_n\}
$ is continuous, then $\{X_n+Y_n\}$ is relatively compact.
\end{corollary}

The following lemma is Problem 3.22(c) in \cite{EK}.
\begin{lemma}\label{L:conlim3}
For fixed $d\ge2$, $\{(X_n^1,\ldots,X_n^d)\}$ is relatively compact in
$D_ {\mathbb{R}^d}[0,\infty)$ if and only if $\{X_n^k\}$ and $\{
X_n^k+X_n^\ell
\}$ are relatively compact in $D_\mathbb{R}[0,\infty)$ for all $k$
and $\ell$.
\end{lemma}

We will also need the following lemma, which connects relative
compactness and convergence in probability. This is Lemma A2.1 in \cite{DK}.
\begin{lemma}\label{L:rcprob}
Let $\{X_n\},X$ be processes with sample paths in
$D_{\mathbb{R}^d}[0,\infty)$ defined on the same probability space.
Suppose that $\{X_n\}$ is relatively compact in $D_{\mathbb{R}^d}
[0,\infty)$ and that for a dense set $H\subset[0,\infty)$,
$X_n(t)\to X(t)$ in probability for all $t\in H$. Then $X_n\to X$
in probability in $D_{\mathbb{R}^d} [0,\infty)$. In particular, if
$X$ is
continuous, then $X_n\to X$ ucp.
\end{lemma}

Our primary tool is the following theorem, which is a special case of
Theorem~2.2 in \cite{KP}.
\begin{theorem}\label{T:KP}
For each $n$, let $Y_n$ be a cadlag, $\mathbb{R}^m$-valued
semimartingale with
respect to a filtration $\{\mathcal{F}_t^n\}$. Suppose that $Y_n=M_n+A_n$,
where $M_n$ is an $\{\mathcal{F}_t^n\}$-local martingale and $A_n$ is
a finite
variation process, and that
%
\begin{equation}\label{KPcond}
\sup_nE\bigl[[M_n]_t+V_t(A_n)\bigr]<\infty
\end{equation}
for each $t\ge0$, where $V_t(A_n)$ is the total variation of $A_n$ on
$[0,t]$ and $[M_n]$ is the quadratic variation of $M_n$. Let $X_n$ be a
cadlag, $\{\mathcal{F}_t^n\}$-adapted, $\mathbb{R}^{k\times
m}$-valued process and define
\[
Z_n(t) = \int_0^t X_n(s-) \,dY_n(s).
\]
Suppose that $(X_n,Y_n)\to(X,Y)$ in law in $D_{\mathbb{R}^{k\times
m}\times
\mathbb{R}
^m} [0,\infty)$. Then, $Y$ is a semimartingale with respect to a
filtration to which $X$ and $Y$ are adapted and $(X_n,Y_n,Z_n)\to
(X,Y,Z)$ in law in $D_ {\mathbb{R}^{k\times m}\times\mathbb
{R}^m\times\mathbb{R}
^k}[0,\infty
)$, where
\[
Z(t) = \int_0^t X(s-) \,dY(s).
\]
If $(X_n,Y_n)\to(X,Y)$ in probability, then $Z_n\to Z$ in probability.
\end{theorem}
\begin{remark}\label{R:KP}
In the setting of Theorem \ref{T:KP}, if $\{W_n\}$ is
another sequence of cadlag, $\{\mathcal{F}_t^n\}$-adapted,
$\mathbb{R}^\ell$-valued processes and $(W_n,X_n,Y_n)$ converges to $(W,X,Y)$
in law
in $D_{\mathbb{R}^\ell\times\mathbb{R}^{k\times m}\times\mathbb
{R}^m}[0,\infty)$, then
$(W_n, X_n,Y_n,Z_n)$ converges to $(W,X,Y,Z)$ in law in
$D_{\mathbb{R}^\ell\times\mathbb{R}^{k\times
m}\times\mathbb{R}^m\times\mathbb{R}^k}[0,\infty)$. This can be
seen by applying
Theorem \ref{T:KP} to $(\overline X_n,\overline Y_n)$, where
$\overline X_n$ is the
block diagonal $(k+\ell)\times(m+1)$ matrix with upper-left entry
$W_n$ and lower-right entry $X_n$, and $\overline Y_n = (0,Y_n^T)^T$.
\end{remark}

\subsection{Estimates from the prequel}\label{sec:kb1}

We now recall some of the basic estimates from \cite{S}.

By (2.6) in \cite{S}, for all $s\leq t$,
\[
\bigl|E|F(t) - F(s)|^2 - (2/\pi)^{1/2}|t - s|^{1/2}\bigr|
\le\pi^{-1/2}(1 + 2^{1/2})^{-1}t^{-3/2}|t - s|^2.
\]
Hence,
%
\begin{equation}\label{sig1}
\pi^{-1/2}|t - s|^{1/2} \le E|F(t) - F(s)|^2
\le2|t - s|^{1/2}.
\end{equation}
In particular, if $\sigma_j^2 = E\Delta F_j^2$, then
%
\begin{equation}\label{sig2}
|\sigma_j^2 - (2/\pi)^{1/2}\Delta t^{1/2}| \le t_j^{-3/2}\Delta t^2
=j^{-3/2}\Delta t^{1/2}
\end{equation}
and
%
\begin{equation}\label{sig3}
\pi^{-1/2}\Delta t^{1/2} \le\sigma_j^2 \le2\Delta t^{1/2}.
\end{equation}
Theorem 2.3 in \cite{S} shows that $F$ has a nontrivial quartic
variation. A special case of this theorem is the fact that
$\sum_{j=1}^{\lfloor nt\rfloor}\Delta F_j^4 \to6t/\pi$ ucp. The
proof can be
easily adapted to show that
%
\begin{equation}\label{quartvar}
\mathop{\sum_{j=1}}_{j\ \mathrm{odd}}^{\lfloor nt\rfloor}\Delta F_{j}^4
\to\frac3\pi t \quad\mbox{and}\quad
\mathop{\sum_{j=1}}_{j\ \mathrm{even}}^{\lfloor nt\rfloor}\Delta F_{j}^4
\to\frac3\pi t
\end{equation}
ucp.

Let
%
\begin{equation}\label{gam}
\gamma_j = 2j^{1/2} - (j-1)^{1/2} - (j+1)^{1/2}
\end{equation}
and note that $\sum_{j=1}^\infty\gamma_j=1$. By (2.4) in \cite{S},
if $i
< j$, then
\[
|E[\Delta F_i\Delta F_j] + (2\pi)^{-1/2}\gamma_{j-i}\Delta t^{1/2}|
\le(t_i + t_j)^{-3/2}\Delta t^2 = (i + j)^{-3/2}\Delta t^{1/2}.
\]
Some related estimates are $0 < \gamma_j \le2^{-1/2}j^{-3/2}$, which
is (2.8) in \cite{S}, and
%
\begin{equation}\label{cross}
-2(t_j - t_i)^{-3/2}\Delta t^2 = -2(j - i)^{-3/2}\Delta t^{1/2}
\le E[\Delta F_i\Delta F_j] < 0,
\end{equation}
which precedes (2.10) in \cite{S}.

Let $\widehat\sigma_j=E[F(t_{j-1})\Delta F_j]$. Since
\begin{eqnarray*}
\widehat\sigma_j + (2\pi)^{-1/2}\Delta t^{1/2} &=& \sum_{i=1}^{j-1}
\bigl(E[\Delta F_i\Delta F_j] + (2\pi)^{-1/2}\gamma_{j-i}\Delta
t^{1/2}\bigr)\\
&&{} + (2\pi)^{-1/2}\Delta t^{1/2}\sum_{i=j}^\infty\gamma_i,
\end{eqnarray*}
it follows that
%
\begin{equation}\label{sighat}
|\widehat\sigma_j + (2\pi)^{-1/2}\Delta t^{1/2}| \le
Cj^{-1/2}\Delta t^{1/2}.
\end{equation}
In particular, $|\widehat\sigma_j|\le C\Delta t^{1/2}$ and
$|\widehat\sigma_j^2 -
(2\pi
)^{-1} \Delta t|\le Cj^{-1/2}\Delta t$.
\begin{lemma}\label{L:sigdel}
If integers $c,i$ and $j$ satisfy $0\le c<i\le j$, then:
\begin{longlist}
\item$|E[(F(t_{i-1}) - F(t_c))\Delta F_j]|
\le C\Delta t^{1/2} ((j - i)\vee1)^{-1/2}$;
\item$|E[(F(t_{j-1}) - F(t_c))\Delta F_i]|
\le C\Delta t^{1/2} [((j - i)\vee1)^{-1/2}+(i-c)^{-1/2}]$;
\item
$|E[F(t_{j-1})\Delta F_i]|\le C\Delta t^{1/2} ((j - i)\vee1)^{-1/2}$.
\end{longlist}
\end{lemma}
\begin{pf}
By (\ref{cross}),
\[
\bigl|E\bigl[\bigl(F(t_{i-1}) - F(t_c)\bigr)\Delta F_j\bigr]\bigr|
\le\sum_{k=c+1}^{i-1} |E[\Delta F_k\Delta F_j]|
\le C\Delta t^{1/2}\sum_{k=c+1}^{i-1}(j - k)^{-3/2}.
\]
Hence,
\[
\bigl|E\bigl[\bigl(F(t_{i-1}) - F(t_c)\bigr)\Delta F_j\bigr]\bigr|
\le C\Delta t^{1/2}\sum_{k=j-i+1}^\infty k^{-3/2},
\]
which proves the first claim.

For the second and third claims, it is easy to see that they hold when
$i\ge j-1$. Assume $i<j-1$. Note that
\begin{eqnarray*}
E[F(t_{j-1})\Delta F_i] &=& \rho(t_i,t_{j-1}) - \rho
(t_{i-1},t_{j-1})\\
&=& \rho(t_{i-1}+\Delta t,t_{j-1}) - \rho(t_{i-1},t_{j-1})\\
&=& \Delta t\partial_s\rho(t_{i-1}+\theta\Delta t,t_{j-1})
\end{eqnarray*}
for some $\theta\in(0,1)$. Since $j>i$, $t_{i-1}+\theta\Delta
t<t_{j-1}$. In the
regime $s<t$, $\partial_s\rho(s,t) = (8\pi)^{-1/2}((t+s)^{-1/2} +
(t-s)^{-1/2})$. Hence, $0<\partial_s \rho(s,t)\le C(t-s)^{-1/2}$. It
follows that
\[
0 < E[F(t_{j-1})\Delta F_i] \le C\Delta t|t_{j-1} - t_i|^{-1/2}
= C\Delta t^{1/2}(j - i - 1)^{-1/2}.
\]
Since $j-i\ge2$, this implies that $E[F(t_{j-1})\Delta F_i] \le
C\Delta
t^{1/2} (j - i)^{-1/2}$, which proves the third claim. Combining this
with the first claim gives
\begin{eqnarray*}
&&
\bigl|E\bigl[\bigl(F(t_{j-1}) - F(t_c)\bigr)\Delta F_i\bigr]\bigr| \\
&&\qquad\le|E[F(t_{j-1})\Delta F_i]|
+ |E[F(t_c)\Delta F_i]|\\
&&\qquad\le C\Delta t^{1/2}[(j - i)^{-1/2} + (i - c)^{-1/2}],
\end{eqnarray*}
which proves the second claim.
\end{pf}

Recall $\gamma_j$, defined by (\ref{gam}). Let
%
\begin{equation}\label{kap}
\kappa= \Biggl(\frac4\pi
+ \frac2\pi\sum_{j=1}^\infty\gamma_j^2(-1)^j \Biggr)^{1/2} > 0
\end{equation}
(the quantity in the brackets is strictly positive by Proposition
4.7 of \cite{S}) and define
%
\begin{equation}\label{Bdef}
B_n(t) = \kappa^{-1}\sum_{j=1}^{2\lfloor nt/2\rfloor} \Delta F_j^2 (-1)^j.
\end{equation}
(Note that this is simply $\kappa^{-1}Q_n^F$, in the notation of Section
\ref{S:intro}.) By Propositions 3.5 and 4.7 in \cite{S},
%
\begin{equation}\label{Bmom}
E|B_n(t) - B_n(s)|^4
\le C \biggl(\frac{2\lfloor nt/2\rfloor- 2\lfloor ns/2\rfloor
}{n} \biggr)^2
\end{equation}
for all $s$ and $t$. Recall that $F(t)=u(x,t)$, where $u$ is given by
(\ref{SPDE}). Let $m$ denote Lebesgue measure and define the filtration
%
\begin{equation}\label{filt}
\mathcal{F}_t = \sigma\{W(A)\dvtx A\subset\mathbb{R}\times
[0,t],m(A)<\infty\}.
\end{equation}
Fix $\tau\ge0$ and define $G(t)=F(t+\tau)-E[F(t+\tau) \mid
\mathcal{F}_\tau]$. In the proof of Lem\-ma~3.6 in \cite{S}, it was shown
that $G$ and $F$ have the same law and that $G$ is independent of
$\mathcal{F}_\tau$. In particular, if $j>c$ and $\Delta\overline
F_j=\Delta F_j-E[\Delta
F_j \mid\mathcal{F}_{t_c}]$, then $\Delta\overline F_j$ is
independent of
$\mathcal{F}_{t_c}$ and equal in law to $\Delta F_{j-c}$.

According to the equation displayed above (3.32) in \cite{S}, if
$0\le\tau\le s\le t$, then
%
\begin{equation}\label{sigcond1}
E|E[F(t)-F(s) \mid\mathcal{F}_\tau]|^2 \le2|t-s|^2|t-\tau|^{-3/2}.
\end{equation}
In particular, $E|\Delta F_j-\Delta\overline F_j|^2\le2\Delta
t^2(t_j-t_c)^{-3/2}=2\Delta
t^ {1/2}(j-c)^{-3/2}$, which, together with (\ref{sig3}) and H\"older's
inequality, implies
that
%
\begin{eqnarray}\label{sigcond2}
E|\Delta F_j^2 - \Delta\overline F{}^2_j|^k
&=& E[|\Delta F_j + \Delta\overline F_j|^k|\Delta F_j - \Delta
\overline F_j|^k]\nonumber\\
&\le& C_k \Delta t^{5k/4} (t_j - t_c)^{-3k/4}\\
&=& C_k \Delta t^{k/2} (j - c)^{-3k/4}.\nonumber
\end{eqnarray}
Finally, we recall the main result of interest to us, which is
Proposition 4.7 in \cite{S}.
\begin{theorem}\label{T:premain}
Let $\{B_n\}$ be given by (\ref{Bdef}) and let $B$ be a standard
Brownian motion, independent of $F$. Then, $(F,B_n)\to(F,B)$ in law in
$D_{\mathbb{R}^2}[0, \infty)$.
\end{theorem}

\subsection{Tools for Gaussian random variables}\label{S:Gaussprelim}

Let
%
\begin{equation}\label{hermite}
h_n(x) = (-1)^n e^{x^2/2}\,\frac{d^n}{dx^n}(e^{-x^2/2})
\end{equation}
be the $n$th Hermite polynomial so that $\{h_n\}$ is an orthogonal
basis of $L^2(\mu)$, where $\mu(dx)=(2\pi)^{-1/2}e^{-x^2/2} dx$; see
Section 1.1.1 of \cite{Nu} for details. Let $\|\cdot\|$ and $\langle
\cdot,\cdot\rangle$ denote the norm and inner product,
respectively, in $L^2(\mu)$.

The first few Hermite polynomials are $h_0(x)=1$, $h_1(x)=x$,
$h_2(x)=x^2-1$ and $h_3(x)=x^3-3x$. We adopt the convention that
$h_{-1}(x)=0$. The Hermite polynomials satisfy the following identities
for $n\ge0$:
%
\begin{eqnarray}
\label{hrm1}
h_n'(x) &=& nh_{n-1}(x),\\
\label{hrm2}
xh_n(x) &=& h_{n+1}(x) + nh_{n-1}(x),\\
\label{hrm3}
h_n(-x) &=& (-1)^nh_n(x).
\end{eqnarray}
Any polynomial can be written as a linear combination of Hermite
polynomials by using the formula
%
\begin{equation}\label{hrmx^n}
x^n = \sum_{j=0}^{\lfloor n/2\rfloor}\pmatrix{n\cr2j}(2j-1)!!h_{n-2j}(x),
\end{equation}
where $(2j-1)!!=(2j-1)(2j-3)(2j-5)\cdots1$. Note that this can be
rewritten as
%
\begin{equation}\label{hrmx^n2}
x^n = \sum_{j=0}^n\pmatrix{n\cr j}E[Y^j]h_{n-j}(x),
\end{equation}
where $Y$ is a standard normal random variable.

In the remaining part of Section \ref{S:Gaussprelim}, $X$ shall
denote a standard normal random variable. If $r\in[-1,1]$, then
$X_r$, $Y_r$ shall denote jointly normal random variables with
mean zero, variance one and $E[X_rY_r]=r$. By Lemma 1.1.1 in
\cite{Nu},
%
\begin{equation}\label{hrmNu}
E[h_n(X_r)h_m(Y_r)] =
\cases{
0, &\quad if $n\ne m$,\cr
n!r^n, &\quad if $n=m$.}
\end{equation}
In particular, $\|h_n\|^2=E[h_n(X)^2]=n!$. Hence, if $g\in L^2(\mu)$, then
%
\begin{equation}\label{hrmex}
g = \sum_{n=0}^\infty\frac1{n!}\langle g,h_n\rangle h_n,
\end{equation}
where the convergence is in $L^2(\mu)$.

If $g$ and $g'$ have polynomial growth and $n\ge1$, then integration by
parts gives
%
\begin{eqnarray}\label{hrmred}\hspace*{34pt}
\langle g,h_n\rangle &=& \frac1{\sqrt{2\pi}}\int
g(x)h_n(x)e^{-x^2/2} \,dx
= \frac{(-1)^n}{\sqrt{2\pi}}\int g(x)\,\frac{d^n}{dx^n}(e^{-x^2/2})\,
dx\nonumber\\[-8pt]\\[-8pt]
&=& \frac{(-1)^{n-1}}{\sqrt{2\pi}}
\int g'(x)\,\frac{d^{n-1}}{dx^{n-1}}(e^{-x^2/2}) \,dx
= \langle g',h_{n-1}\rangle.\nonumber
\end{eqnarray}
That is, $E[g(X)h_n(X)]=E[g'(X)h_{n-1}(X)]$. Using (\ref{hrmNu}) and
(\ref{hrmex}), we can generalize this as follows:
%
\begin{eqnarray}\label{hrmred2}
E[g(X_r)h_n(Y_r)]
&=& \sum_{m=0}^\infty\frac1{m!}\langle g,h_m\rangle
E[h_m(X_r)h_n(Y_r)]\nonumber\\
&=& \langle g,h_n\rangle r^n
= r\langle g',h_{n-1}\rangle r^{n-1}\\
&=& rE[g'(X_r)h_{n-1}(Y_r)].\nonumber
\end{eqnarray}
The following two lemmas will be useful in Section \ref{S:rc}.
\begin{lemma}\label{L:cool}
Suppose $g,h,g',h'$ all have polynomial growth. If $f(r) =
E[g(X_r)h(Y_r)]$, then $f'(r) = E[g'(X_r)h'(Y_r)]$ for all $r\in(-1,1)$.
\end{lemma}
\begin{pf}
By (\ref{hrmex}) and (\ref{hrmNu}), $f(r) = \sum_{n=0}^\infty
\frac
1{n!} \langle g,h_n\rangle\langle h,h_n\rangle r^n$, which, by (\ref
{hrmred}), gives
\begin{eqnarray*}
f'(r) &=& \sum_{n=1}^\infty
\frac1{(n-1)!}\langle g,h _n\rangle\langle h,h_n\rangle r^{n-1}\\
&=& \sum_{n=1}^\infty
\frac1{(n-1)!}\langle g',h_{n-1}\rangle\langle h',h_{n-1}\rangle
r^{n-1}\\
&=& E[g'(X_r)h'(X_r)].
\end{eqnarray*}
\upqed\end{pf}
\begin{lemma}\label{L:cool?}
Suppose $g,g',g'',h,h',h''$ have polynomial growth. Let $U=aX_r$ and
$V=bY_r$. If $\varphi(a,b,r) = E[g(U) h(V)]$, then
\[
\frac{\partial\varphi}{\partial a}(a,b,r)
= aE[g''(U)h(V)] + brE[g'(U)h'(V)]
\]
for all real $a,b$ and all $r\in(-1,1)$.
\end{lemma}
\begin{pf}
By (\ref{hrmex}) and (\ref{hrmNu}), $\varphi(a,b,r) =
\sum_{n=0}^\infty\frac1{n!}\langle g(a\cdot),h_n\rangle\langle
h(b\cdot),h_n\rangle r^n$.
Fix $a_0\in\mathbb{R}$. To justify differentiating under the
summation at
$a_0$, we must show that there exists an $\varepsilon>0$ and a sequence
$C_n(b,r)$ such that
\[
\biggl|\frac{\partial}{\partial a} \biggl[
\frac1{n!}\langle g(a\cdot),h_n\rangle\langle h(b\cdot),h_n\rangle
r^n \biggr] \biggr|
\le C_n(b,r)
\]
for all $|a-a_0|<\varepsilon$, and $\sum_{n=0}^\infty
C_n(b,r)<\infty$. For
this, we use (\ref{hrm2}) and (\ref{hrmred}) to compute
\begin{eqnarray*}
&&\frac{\partial}{\partial a} \biggl[
\frac1{n!}\langle g(a\cdot),h_n\rangle\langle h(b\cdot),h_n\rangle
r^n \biggr]\\
&&\qquad= \frac1{n!}\langle g'(a\cdot),h_{n+1}\rangle\langle h(b\cdot
),h_n\rangle r^n\\
&&\qquad\quad{} + \frac1{(n-1)!}\langle g'(a\cdot),h_{n-1}\rangle\langle
h(b\cdot),h_n\rangle r^n\\
&&\qquad= \frac a{n!}\langle g''(a\cdot),h_n\rangle\langle h(b\cdot
),h_n\rangle r^n\\
&&\qquad\quad{} + \frac b{(n-1)!}
\langle g'(a\cdot),h_{n-1}\rangle\langle h'(b\cdot),h_{n-1}\rangle r^n.
\end{eqnarray*}
Since $|\langle\cdot,h_n/\sqrt{n!}\rangle|\le\|\cdot\|$, we may take
$C_n(b,r)=Mr^n$ for an appropriately chosen constant $M$, provided that
$|r|<1$. We may therefore differentiate under the summation at $a_0$.
Since $a_0$ was arbitrary, we have
\begin{eqnarray*}
\frac{\partial\varphi}{\partial a}(a,b,r) &=& a\sum_{n=0}^\infty
\frac1{n!}\langle g''(a\cdot),h_n\rangle\langle h(b\cdot
),h_n\rangle r^n\\
&&{} + b\sum_{n=1}^\infty\frac1{(n-1)!}
\langle g'(a\cdot),h_{n-1}\rangle\langle h'(b\cdot),h_{n-1}\rangle
r^n\\
&=& aE[g''(U)h(V)] + brE[g'(U)h'(V)]
\end{eqnarray*}
for all $a,b,r$ with $|r|<1$.
\end{pf}

\subsection{Multi-indices and Taylor's theorem}

We recall here the standard multi-index notation. A multi-index is a
vector $\alpha\in\mathbb{Z}_+^d$, where $\mathbb{Z}_+=\mathbb
{N}\cup\{0\}$. We use $e^j$ to
denote the multi-index with $e^j_j =1$ and $e^j_i=0$ for $i\ne j$. If
$\alpha\in\mathbb{Z}_+^d$ and $x\in\mathbb{R}^d$, then
\begin{eqnarray*}
|\alpha| &=& \sum_{j=1}^d \alpha_j,\qquad
\alpha! = \prod_{j=1}^d \alpha_j!,\\
\partial_j &=& \frac{\partial}{\partial x_j},\qquad
\partial^\alpha= \partial_1^{\alpha_1}\cdots
\partial_d^{\alpha_d},\qquad
x^\alpha= \prod_{j=1}^d x_j^{\alpha_j}.
\end{eqnarray*}
Note that by convention, $0^0=1$. Also note that $|x^\alpha
|=y^\alpha$, where
$y_j =|x_j|$ for all~$j$.

Taylor's theorem with integral remainder states that if $g\in
C^{k+1}(\mathbb{R})$, then
%
\begin{equation}\label{Taylor}
g(b) = \sum_{j=0}^k g^{(j)}(a)\frac{(b - a)^j}{j!}
+ \frac1{k!}\int_a^b (b - u)^k g^{(k+1)}(u) \,du.
\end{equation}
Taylor's theorem in higher dimensions is the following.
\begin{theorem}\label{T:Taylor}
If $g\in C^{k+1}(\mathbb{R}^d)$, then
\[
g(b) = \sum_{|\alpha|\le k} \partial^\alpha g(a)\frac{(b -
a)^\alpha}{\alpha!}
+ R,
\]
where
\[
R = (k + 1)\sum_{|\alpha|=k+1}\frac{(b - a)^\alpha}{\alpha!}
\int_0^1 (1 - u)^k \partial^\alpha g\bigl(a + u(b - a)\bigr) \,du.
\]
In particular,
\[
|R| \le(k+1)\sum_{|\alpha|=k+1}M_\alpha|(b-a)^\alpha|,
\]
where $M_\alpha=\sup\{|\partial^\alpha g(a+u(b-a))|\dvtx 0\le u\le1\}$.
\end{theorem}

For integers $a$ and $b$ with $a\ge0$, we adopt the convention that
\[
\pmatrix{a\cr b} =
\cases{
{\dfrac{a!}{b!(a-b)!}},&\quad if $0\le b\le a$,\vspace*{2pt}\cr
0,&\quad if $b<0$ or $b>a$.}
\]
We define
\[
\pmatrix{\gamma\cr\alpha} = \prod_{j=1}^d\pmatrix{\gamma_j\cr
\alpha_j}
\]
for any multi-indices $\gamma$ and $\alpha$. Later in the paper, we shall
need the following two combinatorial lemmas.
\begin{lemma}\label{L:combo1}
Let $a,b$ and $c$ be integers. If $a\ge0$ and $0\le c\le a$,
then
\[
\sum_{j=0}^c\pmatrix{a-c\cr b-j}\pmatrix{c\cr j} = \pmatrix{a\cr b}.
\]
\end{lemma}
\begin{pf}
The proof is by induction on $a$. For $a=0$, the lemma is trivial.
Suppose the lemma holds for $a-1$. Since the lemma clearly holds for
$c=0$ or $c=a$, we may assume $0<c\le a-1$. In that case,
\begin{eqnarray*}
\pmatrix{a\cr b} &=& \pmatrix{a-1\cr b} + \pmatrix{a-1\cr b-1}\\
&=& \sum_{j=0}^c \left[\pmatrix{a-1-c\cr b-j}
+ \pmatrix{a-1-c\cr b-1-j} \right]\pmatrix{c\cr j}\\
&=& \sum_{j=0}^c\pmatrix{a-c\cr b-j}\pmatrix{c\cr j}.\mbox{\hspace*{150pt}\qed\hspace*{-150pt}}
\end{eqnarray*}
\noqed\end{pf}

Suppose $\alpha$ and $\gamma$ are multi-indices. We will write
$\alpha\le\gamma$ if $\alpha_j\le\gamma_j$ for all $j$.
\begin{lemma}\label{L:combo2}
If $\gamma$ is a multi-index in $\mathbb{Z}_+^d$ and $m\ge0$, then
\[
\mathop{\sum_{|\alpha|=m}}_{\alpha\le\gamma}\pmatrix{\gamma
\cr\alpha}
= \pmatrix{|\gamma|\cr m}.
\]
\end{lemma}
\begin{pf}
We shall prove this by induction on $d$. If $d=1$, then the lemma
is trivial. Suppose the lemma is true for $d-1$. Let $\gamma$ be a
multi-index in $\mathbb{Z}_+^d$ and fix $m$ with $0\le m\le|\gamma
|$. For
multi-indices $\alpha$ and $\gamma$, let
$\widehat\alpha=(\alpha_1,\ldots,\alpha_{d-1})$ and
$\widehat\gamma=(\gamma_1,\ldots,\gamma_{d-1})$. Then,
\begin{eqnarray*}
\mathop{\sum_{|\alpha|=m}}_{\alpha\le\gamma}\pmatrix{\gamma
\cr\alpha}
&=& \sum_{\alpha_d=0}^{m\wedge\gamma_d}
\mathop{\sum_{|\widehat\alpha|=m-\alpha_d}}_{\widehat\alpha
\le\widehat\gamma}
\pmatrix{\widehat\gamma\cr\widehat\alpha}\pmatrix{\gamma
_d\cr\alpha_d}\\
&=& \sum_{\alpha_d=0}^{m\wedge\gamma_d}
\pmatrix{|\widehat\gamma|\cr m-\alpha_d}\pmatrix{\gamma_d\cr\alpha_d}\\
&=& \sum_{\alpha_d=0}^{\gamma_d}
\pmatrix{|\gamma|-\gamma_d\cr m-\alpha_d}\pmatrix{\gamma_d\cr
\alpha_d}.
\end{eqnarray*}
Applying Lemma \ref{L:combo1} completes the proof.
\end{pf}

\section{Fourth order integrals}\label{S:4ints}

\begin{theorem}\label{T:4ints}
Suppose $g\dvtx\mathbb{R}\times[0,\infty)\to\mathbb{R}$ is continuous.
For each $n$, let
$\{s_j^*\}$ and $\{t_j^*\}$ be collections of points with
$s_j^*,t_j^*\in[t_{j-1},t_j]$. Then,
%
\begin{eqnarray}\label{4ints}
\lim_{n\to\infty}\mathop{\sum_{j=1}}_{j\ \mathrm{odd}}^{\lfloor
nt\rfloor}
g(F(s_j^*),t_j^*)\Delta F_j^4
&=& \lim_{n\to\infty}\mathop{\sum_{j=1}}_{j\ \mathrm{even}}^{\lfloor
nt\rfloor}
g(F(s_j^*),t_j^*)\Delta F_j^4\nonumber\\[-8pt]\\[-8pt]
&=& \frac3\pi\int_0^t g(F(s),s) \,ds,\nonumber
\end{eqnarray}
where the convergence is ucp.
\end{theorem}
\begin{pf}
We prove only the first limit. The proof for the other limit
is nearly identical. Let
\[
X_n(t) = \sum_{j=1}^\infty g(F(s_j^*),t_j^*)1_{[t_{j-1},t_j)}(t)
\]
and
\[
A_n(t) = \mathop{\sum_{j=1}}_{j\ \mathrm{odd}}^{\lfloor nt\rfloor}
\Delta F_j^4
\]
so that
\[
\mathop{\sum_{j=1}}_{j\ \mathrm{odd}}^{\lfloor nt\rfloor}
g(F(s_j^*),t_j^*)\Delta F_j^4 = \int_0^t X_n(s-) \, dA_n(s).
\]
By (\ref{quartvar}), $A_n(t)\to3t/\pi$ ucp. Also, by the continuity of
$g$ and $F$, $X_n\to g(F(\cdot),\cdot)$ ucp. Finally, note that the
expected total variation $V_t(A_n)$ of $A_n$ on $[0,t]$ is uniformly
bounded in $n$. That is,
\[
E[V_t(A_n)] = \mathop{\sum_{j=1}}_{j\ \mathrm{odd}}^{\lfloor nt\rfloor}
E\Delta F_j^4 \le C\sum_{j=1}^{\lfloor nt\rfloor}\Delta t \le Ct.
\]
By Theorem \ref{T:KP}, (\ref{4ints}) holds with the convergence
being in probability in $D_\mathbb{R}[0,\infty)$. Since the limit is
continuous, (\ref{4ints}) holds ucp.
\end{pf}

If $r$ and $k$ are nonnegative integers with $r\le k$, then we shall
use the notation $g\in C_r^{k,1}(\mathbb{R}\times[0,\infty))$ to
mean that
%
\begin{eqnarray}
\label{c1k0}\quad
&&g\dvtx\mathbb{R}\times[0,\infty)\to\mathbb{R}\mbox{ is
continuous,}\\
\label{c1k1}
&&\partial_x^j g \mbox{ exists and is continuous on $\mathbb
{R}\times[0,\infty)$}\qquad
\mbox{for all $0\le j\le k$},\\
\label{c1k2}
&&\partial_t\partial_x^j g \mbox{ exists and is continuous on
$\mathbb{R}\times
(0,\infty)$}\qquad\mbox{for all $0\le j\le r$},\\
\label{c1k3}
&&\overline{\lim_{t\to0}} \sup_{x\in K}|\partial_t\partial_x^j
g(x,t)| \,dt<\infty\nonumber\\[-8pt]\\[-8pt]
\eqntext{\mbox{for all compact $K \subset\mathbb{R}$ and all $0\le j\le
r$}.}
\end{eqnarray}
Note that $g\in C_r^{k,1}$ implies $\partial_x^j g\in C_{r-j}^{k-j,1}$
whenever $r\ge j$. For functions of one spatial dimension, we shall
henceforth use standard prime notation to denote spatial derivatives.
For example, $g''=\partial_x^2g$ and $g^{(4)}=\partial_x^4g$.

Typically, we shall need (\ref{c1k2}) and (\ref{c1k3}) only when $j=0$.
There are a few places, however, where $j>0$ is needed. We need $j=3$
in the derivation of (\ref{kb5.2}), which is used in the proofs of both
Theorem \ref{T:trapexpan} and Corollary \ref{C:expan2}; we need $j=2$
in the proof of Lemma \ref{L:claim02}; we need $j=4$ in the proof of
Theorem \ref{T:main}.
Note that $\partial_t\partial_x^j g$ need not be continuous at
$t=0$. In
particular, $\partial_t\partial_x^j g$ need not be bounded on sets
of the form
$K\times(0,\varepsilon]$.

Recall that $X_n(t)\approx Y_n(t)$ means that $X_n-Y_n\to0$ ucp.
\begin{theorem}\label{T:expan1}
If $g\in C^{5,1}_0(\mathbb{R}\times[0,\infty))$, then
\begin{eqnarray*}
I_n(g',t) &\approx& g(F(t),t) - g(F(0),0) - \int_0^t \partial
_tg(F(s),s)
\,ds\\
&&{} - \frac12\sum_{j=1}^{\lfloor nt/2\rfloor}
g''(F(t_{2j-1}),t_{2j-1})(\Delta F_{2j}^2 - \Delta F_{2j-1}^2)\\
&&{} - \frac16\sum_{j=1}^{\lfloor nt/2\rfloor}g'''(F(t_{2j-1}),t_{2j-1})
(\Delta F_{2j}^3 + \Delta F_{2j-1}^3),
\end{eqnarray*}
where $I_n(g,t)$ is given by (\ref{riem}).
\end{theorem}
\begin{pf}
By (\ref{Taylor}),
\begin{eqnarray*}
&&
g(x + h_1,t) - g(x + h_2,t)\\
&&\qquad= \sum_{j=1}^4\frac1{j!} g^{(j)}(x,t)(h_1^j - h_2^j)
+ R(x,h_1,t) - R(x,h_2,t),
\end{eqnarray*}
where
\[
R(x,h,t) = \frac1{4!} \int_0^h{(h-u)^4g^{(5)}(x + u,t) \,du}.
\]
Taking $x=F(t_{2j-1})$, $h_1=\Delta F_{2j}$ and $h_2=-\Delta
F_{2j-1}$, we have
\begin{eqnarray*}
&& g(F(t_{2j}),t_{2j-1}) - g(F(t_{2j-2}),t_{2j-1})\\
&&\qquad = \sum_{j=1}^4\frac1{j!} g^{(j)}(F(t_{2j-1}),t_{2j-1})
\bigl(\Delta F_{2j}^j - (-1)^j\Delta F_{2j-1}^j\bigr)\\
&&\qquad\quad{} + R(F(t_{2j-1}),\Delta F_{2j},t_{2j-1})\\
&&\qquad\quad{} - R(F(t_{2j-1}),-\Delta F_{2j-1},t_{2j-1}).
\end{eqnarray*}
Let $N(t)=2\lfloor nt/2\rfloor/n$. That is, if $t\in
[t_{2j-2},t_{2j})$, then
$N(t)= t_{2j-2}$. Let $F_n(t)=F(N(t))$. Then,
\begin{eqnarray*}
g(F(t_{2j}),t_{2j}) - g(F(t_{2j}),t_{2j-1})
&=& \int_{t_{2j-1}}^{t_{2j}}\partial_tg\bigl(F_n(s + \Delta t),s\bigr) \,ds,\\
g(F(t_{2j-2}),t_{2j-1}) - g(F(t_{2j-2}),t_{2j-2})
&=& \int_{t_{2j-2}}^{t_{2j-1}}\partial_tg(F_n(s),s) \,ds\\
&=& \int_{t_{2j-2}}^{t_{2j-1}}\partial_tg\bigl(F_n(s + \Delta t),s\bigr)
\,ds.
\end{eqnarray*}
Thus,
\begin{eqnarray*}
g(F(t),t) &=& g(F(0),0) + \sum_{j=1}^{\lfloor nt/2\rfloor}
\{g(F(t_{2j}),t_{2j}) - g(F(t_{2j-2}),t_{2j-2})\}\\
&&{} + g(F(t),t) - g(F_n(t),N(t))\\
&=& g(F(0),0) + \int_0^{N(t)}\partial_tg\bigl(F_n(s + \Delta t),s\bigr) \,ds
+ I_n(g',t)\\
&&{} + \frac12 \sum_{j=1}^{\lfloor nt/2\rfloor}
g''(F(t_{2j-1}),t_{2j-1})(\Delta F_{2j}^2 - \Delta F_{2j-1}^2)\\
&&{} + \frac16 \sum_{j=1}^{\lfloor nt/2\rfloor}
g'''(F(t_{2j-1}),t_{2j-1})(\Delta F_{2j}^3 + \Delta F_{2j-1}^3)\\
&&{} +
\varepsilon_n(g,t),
\end{eqnarray*}
where
%
\begin{eqnarray}\label{kb4}
\varepsilon_n(g,t) &=& \frac1{24}\sum_{j=1}^{\lfloor nt/2\rfloor}
g^{(4)}(F(t_{2j-1}),t_{2j-1})(\Delta F_{2j}^4 - \Delta
F_{2j-1}^4)\nonumber\\
&&{} + \sum_{j=1}^{\lfloor nt/2\rfloor}
\{R(F(t_{2j-1}),\Delta F_{2j},t_{2j-1})\nonumber\\[-8pt]\\[-8pt]
&&\hspace*{38.8pt}{} - R(F(t_{2j-1}),-\Delta F_{2j-1},t_{2j-1})\}\nonumber\\
&&{} + g(F(t),t) - g(F_n(t),N(t)).\nonumber
\end{eqnarray}
By (\ref{c1k2}), (\ref{c1k3}), the continuity of $F$ and dominated
convergence,
\[
\int_0^{N(t)}\partial_tg\bigl(F_n(s + \Delta t),s\bigr) \,ds
\to\int_0^t \partial_t g(F(s),s) \,ds
\]
uniformly on compacts, with probability one. Therefore, it will suffice
to show that $\varepsilon_n(g,t)\to0$ ucp.

First, assume that $g$ has compact support. By the continuity of $g$
and the almost sure continuity of $F$, $g(F(t),t) - g(F_n(t), N(t)) \to
0$ ucp. Since $g^{(5)}$ is bounded, $|R(x,h,t)|\le C|h|^5$. Thus,
\begin{eqnarray*}
&&\Biggl|\sum_{j=1}^{\lfloor nt/2\rfloor}\{R(F(t_{2j-1}),\Delta
F_{2j},t_{2j-1})
- R(F(t_{2j-1}),-\Delta F_{2j-1},t_{2j-1})\} \Biggr|\\
&&\qquad\le C\sum_{j=1}^{\lfloor nt/2\rfloor}|\Delta F_j|^5
\end{eqnarray*}
and
\begin{eqnarray*}
E \Biggl[\sup_{0\le t\le T}\sum_{j=1}^{\lfloor nt/2\rfloor}|\Delta
F_j|^5 \Biggr]
&=& \sum_{j=1}^{\lfloor nT/2\rfloor}E|\Delta F_j|^5
= C\sum_{j=1}^{\lfloor nT/2\rfloor}\sigma_j^5\\
&\le& CnT\Delta t^{5/4} = CT\Delta t^{1/4}.
\end{eqnarray*}
It follows that
\[
\sum_{j=1}^{\lfloor nt/2\rfloor}\bigl\{R(F(t_{2j-1}),\Delta F_{2j},t_{2j-1})
- R(F(t_{2j-1}),-\Delta F_{2j-1},t_{2j-1})\bigr\} \to0
\]
ucp. An application of Theorem \ref{T:4ints} to the first sum in
(\ref{kb4}) completes the proof that $\varepsilon_n(g,t) \to0$ ucp,
in the case where
$g$ has compact support.

To deal with the general case, we use the following truncation
argument, which we will make use of several times throughout this
paper. Fix $T>0$ and $\eta>0$. Choose $L>T$ so large that
\[
P \Bigl({\sup_{0\le t\le T}} |F(t)| \ge L \Bigr) < \eta.
\]
Let $\varphi\in C^\infty(\mathbb{R})$ have compact support with
$\varphi\equiv1$ on
$[-L, L]$. Define $h(x,t)=g(x,t)\varphi(x)\varphi(t)$. Then, $h\in
C^{5,1}_0(\mathbb{R}
\times[0, \infty))$, $h$ has compact support and $h=g$ on
$[-L,L]\times
[0,T]$. By the above, we may choose $n_0$ such that
\[
P \Bigl({\sup_{0\le t\le T}} |\varepsilon_n(h,t)| > \eta\Bigr)
< \eta
\]
for all $n\ge n_0$. Hence,
\begin{eqnarray*}
&&P \Bigl({\sup_{0\le t\le T}} |\varepsilon_n(g,t)| > \eta\Bigr)
\\
&&\qquad\le P \Bigl({\sup_{0\le t\le T}} |F(t)| \ge L \Bigr)
+ P \Bigl({\sup_{0\le t\le T}} |\varepsilon_n(h,t)| > \eta\Bigr)\\
&&\qquad< 2\eta
\end{eqnarray*}
for all $n\ge n_0$, which shows that $\varepsilon_n(g,t)\to0$ ucp
and completes
the proof.
\end{pf}
\begin{theorem}\label{T:trapexpan}
If $g\in C^{5,1}_3(\mathbb{R}\times[0,\infty))$, then
\begin{eqnarray*}
T_n(g',t) &\approx& g(F(t),t) - g(F(0),0)
- \int_0^t \partial_tg(F(s),s) \,ds\\
&&{} + \frac1{24}\sum_{j=1}^{\lfloor nt\rfloor}g'''(F(t_j),t_j)
(\Delta F_{j+1}^3 + \Delta F_j^3),
\end{eqnarray*}
where $T_n(g,t)$ is given by (\ref{Triem}).
\end{theorem}
\begin{pf}
As in the proof of Theorem \ref{T:expan1}, we may assume $g$
has compact support. Define
\[
\widehat I_n(g,t) = \sum_{j=1}^{\lfloor nt/2\rfloor}g(F(t_{2j}),t_{2j})
\bigl(F(t_{2j+1}) - F(t_{2j-1})\bigr).
\]
The proof of Theorem \ref{T:expan1} can be easily adapted to show that
%
\begin{eqnarray}\label{offexpan}
\widehat I_n(g',t) &\approx& g(F(t),t) - g(F(0),0)
- \int_0^t \partial_tg(F(s),s) \,ds\nonumber\\
&&{} - \frac12\sum_{j=1}^{\lfloor nt/2\rfloor}
g''(F(t_{2j}),t_{2j})(\Delta F_{2j+1}^2 - \Delta F_{2j}^2)\\
&&{} - \frac16\sum_{j=1}^{\lfloor nt/2\rfloor}g'''(F(t_{2j}),t_{2j})
(\Delta F_{2j+1}^3 + \Delta F_{2j}^3).\nonumber
\end{eqnarray}
Note that
\begin{eqnarray*}
\hspace*{-4pt}&&
I_n(g',t) + \widehat I_n(g',t) \\
\hspace*{-4pt}&&\qquad= \mathop{\sum_{j=1}}_
{j\ \mathrm{odd}}^{2\lfloor nt/2\rfloor}g'(F(t_j),t_j)(\Delta F_{j+1} +
\Delta F_j) + \mathop{\sum_{j=1}}_{
j\ \mathrm{even}}^{2\lfloor nt/2\rfloor}g'(F(t_j),t_j)(\Delta F_{j+1}
+ \Delta F_j)\\
\hspace*{-4pt}&&\qquad= \sum_{j=1}^{2\lfloor nt/2\rfloor}g'(F(t_j),t_j)(\Delta F_{j+1} +
\Delta F_j).
\end{eqnarray*}
Also, note that
\[
T_n(g',t) = \frac12 \Biggl(
\sum_{j=0}^{\lfloor nt\rfloor-1}g'(F(t_j),t_j)\Delta F_{j+1}
+ \sum_{j=0}^{\lfloor nt\rfloor}g'(F(t_j),t_j)\Delta F_j \Biggr).
\]
By the continuity of $F$ and $g'$, this shows that
\[
T_n(g',t) \approx\frac{I_n(g',t) + \widehat I_n(g',t)}2.
\]
By (\ref{offexpan}) and Theorem \ref{T:expan1},
we have
%
\begin{eqnarray}\label{kb5}
T_n(g',t) &\approx& g(F(t),t) - g(F(0),0)
- \int_0^t \partial_tg(F(s),s) \,ds\nonumber\\
&&{} + \frac14 \sum_{j=1}^{\lfloor nt\rfloor}
\bigl(g''(F(t_j),t_j) - g''(F(t_{j-1}),t_{j-1})\bigr)\Delta F_j^2\\
&&{} - \frac1{12}\sum_{j=1}^{\lfloor nt\rfloor}
g'''(F(t_j),t_j)(\Delta F_{j+1}^3 + \Delta F_j^3).\nonumber
\end{eqnarray}
Since $g\in C^{5,1}_3(\mathbb{R}\times[0,\infty))$, we may use the Taylor
expansion $f(b)-f(a)=\frac12 (f'(a)+f'(b))(b-a)+O(|b-a|^3)$ with
$f=g''$ to obtain
%
\begin{eqnarray}\label{kb5.1}
&&
g''(F(t_j),t_j) - g''(F(t_{j-1}),t_{j-1})
\nonumber\\
&&\qquad= \int_{t_{j-1}}^{t_j}\partial_tg''(F(t_j),s) \,ds\nonumber\\
&&\qquad\quad{} + \frac12\bigl(g'''(F(t_{j-1}),t_{j-1})
+ g'''(F(t_j),t_{j-1})\bigr)\Delta F_j + R\\
&&\qquad= \int_{t_{j-1}}^{t_j}\partial_tg''(F(t_j),s) \,ds
- \frac12\Delta F_j\int_{t_{j-1}}^{t_j}\partial_tg'''(F(t_j),s)
\,ds\nonumber\\
&&\qquad\quad{} + \frac12\bigl(g'''(F(t_{j-1}),t_{j-1})
+ g'''(F(t_j),t_j)\bigr)\Delta F_j + R,\nonumber
\end{eqnarray}
where $|R|\le C|\Delta F_j|^3$. Since $g$ has compact support, we may use
(\ref{c1k3}) with $K=\mathbb{R}$ and $j=3$ to conclude that the above
integrals are bounded by $C \Delta t$. This yields
%
\begin{eqnarray}\label{kb5.2}
&&
\sum_{j=1}^{\lfloor nt\rfloor}
\bigl(g''(F(t_j),t_j) - g''(F(t_{j-1}),t_{j-1})\bigr)\Delta
F_j^2\nonumber\\[-8pt]\\[-8pt]
&&\qquad= \sum_{j=1}^{\lfloor nt\rfloor} \frac12
\bigl(g'''(F(t_j),t_j) + g'''(F(t_{j-1}),t_{j-1})\bigr)\Delta F_j^3
+ \widetilde R,\nonumber
\end{eqnarray}
where $|\widetilde R|\le C\sum(\Delta t\Delta F_j^2 + \Delta
t|\Delta F_j|^3 + |\Delta
F_j|^5)$. We can combine this formula with (\ref{kb5}) to complete the proof.
\end{pf}

\section{Third order integrals}\label{S:3ints}

To analyze the third order integrals, we will need a Taylor
expansion of a different kind. That is, we will need an expansion
for the expectation of functions of jointly Gaussian random
variables. For this Gaussian version of Taylor's theorem, we first
introduce some terminology. We shall say that a function
$g\dvtx\mathbb{R}^d\to\mathbb{R}$ has \textit{polynomial growth} if there exist
positive constants $K$ and $r$ such that
\[
|g(x)| \le K(1 + |x|^r)
\]
for all $x\in\mathbb{R}^d$. If $k$ is nonnegative integer, we shall
say that a
function $g$ has \textit{polynomial growth of order $k$} if $g\in
C^k(\mathbb{R}
^d)$ and there exist positive constants $K$ and $r$ such that
\[
|\partial^\alpha g(x)| \le K(1 + |x|^r)
\]
for all $x\in\mathbb{R}^d$ and all $|\alpha|\le k$.
\begin{theorem}\label{T:GaussTaylor}
Let $k$ be a nonnegative integer. Suppose $h\dvtx\mathbb{R}\to\mathbb
{R}$ is measurable
and has polynomial growth, and $f\in C^{k+1}(\mathbb{R}^d)$ has polynomial
growth of order $k+1$, both with common constants $K$ and $r$. Suppose,
also, that $\partial^\alpha f$ has polynomial growth with constants
$K_\alpha$ and
$r$ for all $|\alpha| \le k+1$. Let $\xi\in\mathbb{R}^d$ and $Y\in
\mathbb{R}$ be jointly
normal with mean zero. Suppose that $EY^2=1$ and $E\xi_j^2\le\nu$ for
some $\nu>0$. Define $\rho\in\mathbb{R}^d$ by $\rho_j=E[\xi_jY]$. Then,
\[
E[f(\xi)h(Y)] = \sum_{|\alpha|\le k}\frac1{\alpha!} \rho
^\alpha
E[\partial^\alpha f(\xi- \rho Y)]E\bigl[Y^{|\alpha|}h(Y)\bigr] + R,
\]
where $|R|\le C\sum_{|\alpha|=k+1}K_\alpha|\rho^\alpha|$ and
$C$ depends
only on
$K$, $r$, $\nu$, $k$ and $d$. In particular, $|R|\le C|\rho|^{k+1}$.
\end{theorem}
\begin{pf}
Let $U=\xi-\rho Y$ and define $\varphi\dvtx\mathbb{R}^d\to\mathbb
{R}$ by $\varphi(x) =
E[f(U + xY)h(Y)]$. Since $h$ and $f$ have polynomial growth and
all derivatives of $f$ up to order $k+1$ have polynomial
growth,\vspace*{1pt}
we may differentiate under the expectation and conclude that
$\varphi\in C^{k+1}(\mathbb{R}^d)$. Hence, by Theorem \ref
{T:Taylor} and the
fact that $U$ and $Y$ are independent,
\begin{eqnarray*}
E[f(\xi)h(Y)] &=& \varphi(\rho) = \sum_{|\alpha|\le k}
\frac1{\alpha!} \rho^\alpha\partial^\alpha\varphi(0) + R\\
&=& \sum_{|\alpha|\le k}\frac1{\alpha!} \rho^\alpha
E[\partial^\alpha f(U)]E\bigl[Y^{|\alpha|}h(Y)\bigr] + R,
\end{eqnarray*}
where
\[
|R| \le(k+1)\sum_{|\alpha|=k+1}M_\alpha|\rho^\alpha|
\]
and $M_\alpha=\sup\{|\partial^\alpha\varphi(u\rho)|\dvtx 0\le u\le
1\}$. Note that
\[
\partial^\alpha\varphi(u\rho) = E[\partial^\alpha f(U + u\rho
Y)Y^{k+1}h(Y)]
= E\bigl[\partial^\alpha f\bigl(\xi- \rho(1 - u)Y\bigr)Y^{k+1}h(Y)\bigr].
\]
Hence,
\begin{eqnarray*}
|\partial^\alpha\varphi(u\rho)|
&\le& K_\alpha K E\bigl[\bigl(1 + |\xi- \rho(1 - u)Y|^r\bigr)|Y|^{k+1}(1 + |Y|^r)\bigr]\\
&\le& K_\alpha K E[(1 + 2^r|\xi|^r + 2^r|\rho|^r|Y|^r)
(|Y|^{k+1} + |Y|^{k+1+r})].
\end{eqnarray*}
Since $|\rho|^2\le\nu d$, this shows that $M_\alpha\le CK_\alpha$ and
completes the proof.
\end{pf}
\begin{corollary}\label{C:GaussTaylor}
Recall the Hermite polynomials $h_n(x)$ from (\ref{hermite}). Under the
hypotheses of Theorem \ref{T:GaussTaylor},
\[
E[f(\xi)h(Y)] = \sum_{|\alpha|\le k}\frac1{\alpha!} \rho
^\alpha
E[\partial^\alpha f(\xi)]E\bigl[h_{|\alpha|}(Y)h(Y)\bigr] + R,
\]
where $|R|\le C\sum_{|\alpha|=k+1}K_\alpha|\rho^\alpha|$ and
$C$ depends
only on
$K$, $r$, $\nu$, $k$ and $d$. In particular, $|R|\le C|\rho|^{k+1}$.
\end{corollary}
\begin{pf}
Recursively define the sequences $\{a^{(n)}_j\}_{j=0}^\infty$
by $a^{(0)}_j =E[Y^jh(Y)]$ and
%
\begin{equation}\label{coeffs}
a^{(n+1)}_j =
\cases{
a^{(n)}_j,&\quad if $j\le n$,\cr
a^{(n)}_j - {\pmatrix{j\cr n}}a^{(n)}_nE[Y^{j-n}],&\quad if $j\ge n+1$.}
\end{equation}
We will show that for all $0\le n\le k+1$,
%
\begin{eqnarray}\label{star01}
E[f(\xi)h(Y)] &=& \sum_{|\alpha|\le n-1}
\frac1{\alpha!}\rho^\alpha E[\partial^\alpha f(\xi
)]a_{|\alpha|}^{(n)}\nonumber\\[-8pt]\\[-8pt]
&&{} + \sum_{n\le|\alpha|\le k}
\frac1{\alpha!}\rho^\alpha E[\partial^\alpha f(\xi- \rho
Y)]a_{|\alpha|}^{(n)}
+ R,\nonumber
\end{eqnarray}
where $|R|\le C\sum_{|\alpha|=k+1}K_\alpha|\rho^\alpha|$ and
$C$ depends
only on
$K$, $r$,
$\nu$, $k$ and $d$. The proof is by induction on $n$. The case
$n=0$ is given by Theorem \ref{T:GaussTaylor}. Suppose
(\ref{star01}) holds for some $n<k+1$. Fix $\alpha$ such that
$|\alpha|=n$. Let $c_k $ denote $ E [Y^k]$. Applying Theorem
\ref{T:GaussTaylor} to $\partial^\alpha f$ with $h(y)=1$ gives
\begin{eqnarray*}
E[\partial^\alpha f(\xi)] &=& \sum_{|\beta|\le k-n}\frac1{\beta
!} \rho^\beta
E[\partial^{\alpha+\beta}f(\xi- \rho Y)]c_{|\beta|} +
\widehat R_\alpha\\
&=& E[\partial^\alpha f(\xi- \rho Y)]
+ \sum_{1\le|\beta|\le k-n}\frac1{\beta!} \rho^\beta
E[\partial^{\alpha+\beta}f(\xi- \rho Y)]c_{|\beta|} +
\widehat R_\alpha,
\end{eqnarray*}
where $|\widehat R_\alpha|\le C\sum_{|\beta|=k+1-n}K_{\alpha
+\beta}|\rho^\beta|$. Hence,
by (\ref{star01}),
%
\begin{eqnarray}\label{star02}
E[f(\xi)h(Y)] &=& \sum_{|\alpha|\le n}
\frac1{\alpha!}\rho^\alpha E[\partial^\alpha f(\xi
)]a_{|\alpha|}^{(n)}\nonumber\\
&&{} + \sum_{n+1\le|\alpha|\le k}
\frac1{\alpha!}\rho^\alpha E[\partial^\alpha f(\xi- \rho
Y)]a_{|\alpha|}^{(n)}\\
&&{} - S + R^*,\nonumber
\end{eqnarray}
where
\begin{eqnarray*}
|R^*| &\le& |R| + C\sum_{|\alpha|=n}|\rho^\alpha||\widehat
R_\alpha|\\
&\le& |R| + C\sum_{|\alpha|=n}|\rho^\alpha|
\sum_{|\beta|=k+1-n} K_{\alpha+\beta}|\rho^\beta|\\
&\le& C\sum_{|\alpha|=n} K_\alpha|\rho^\alpha|
\end{eqnarray*}
and
\[
S = \sum_{|\alpha|=n}\sum_{1\le|\beta|\le k-n}
\frac1{\alpha!\beta!}\rho^{\alpha+\beta}E[\partial^{\alpha
+\beta}f(\xi- \rho Y)]
a_n^{(n)}c_{|\beta|}.
\]
Making the change of index $\gamma=\alpha+\beta$ and using Lemma
\ref
{L:combo2} gives
\begin{eqnarray*}
S &=& \sum_{n+1\le|\gamma|\le k}
\mathop{\sum_{|\alpha|=n}}_{\alpha\le\gamma}
\pmatrix{\gamma\cr\alpha}\frac1{\gamma!}\rho^\gamma E[\partial
^\gamma f(\xi-\rho Y)]
a_n^{(n)}c_{|\gamma|-n}\\
&=& \sum_{n+1\le|\gamma|\le k}\pmatrix{|\gamma|\cr n}
\frac1{\gamma!}\rho^\gamma E[\partial^\gamma f(\xi-\rho Y)]
a_n^{(n)}c_{|\gamma|-n}.
\end{eqnarray*}
Substituting this into (\ref{star02}) and using (\ref{coeffs}) shows that
%
\begin{eqnarray}
E[f(\xi)h(Y)] &=& \sum_{|\alpha|\le n}
\frac1{\alpha!}\rho^\alpha E[\partial^\alpha f(\xi
)]a_{|\alpha|}^{(n+1)}\nonumber\\[-8pt]\\[-8pt]
&&{} + \sum_{n+1\le|\alpha|\le k}
\frac1{\alpha!}\rho^\alpha E[\partial^\alpha f(\xi- \rho
Y)]a_{|\alpha|}^{(n+1)}
+ R^*,\nonumber
\end{eqnarray}
which completes the induction.

By (\ref{star01}) with $n=k+1$, it remains only to show that
%
\begin{equation}\label{newmagic}
a_j^{(n)} = E[h_j(Y)h(Y)] \qquad\mbox{for all $j\le n$}.
\end{equation}
The proof is by induction on $n$. For $n=0$, the claim is trivial.
Suppose (\ref{newmagic}) holds for all $n\le N$. If $j \leq N$,
then (\ref{coeffs}) implies $a_j^{(N+1)}=a_j^{(N)}=E[h_j(Y)h(Y)]$.
If $j=N+1$, then
\[
a_{N+1}^{(N+1)} = a_{N+1}^{(N)} - \pmatrix{N+1\cr N}
a_{N}^{(N)}E[Y].
\]
Using induction, this gives
\begin{eqnarray*}
a_{N+1}^{(N+1)} &=& a_{N+1}^{(0)} - \sum_{j=0}^{N}\pmatrix{N+1\cr j}
a_j^{(j)}E[Y^{N+1-j}]\\
&=& E[Y^{N+1}h(Y)] - \sum_{j=0}^{N} \pmatrix{N+1\cr j}
E[h_j(Y)h(Y)]E[Y^{N+1-j}]\\
&=& E \Biggl[ \Biggl\{Y^{N+1} - \sum_{j=0}^{N}\pmatrix{N+1\cr j}
E[Y^{N+1-j}]h_j(Y) \Biggr\}h(Y) \Biggr].
\end{eqnarray*}
By (\ref{hrmx^n2}),
\[
Y^{N+1} = \sum_{j=0}^{N+1}\pmatrix{N+1\cr j}E[Y^j]h_{N+1-j}(Y)
= \sum_{j=0}^{N+1}\pmatrix{N+1\cr j}E[Y^{N+1-j}]h_j(Y).
\]
Hence, $a_{N+1}^{(N+1)}=E[h_{N+1}(Y)h(Y)]$, completing the proof
of (\ref{newmagic}).
\end{pf}
\begin{theorem}\label{T:3ints}
If $g\in C^{4,1}_0(\mathbb{R}\times[0,\infty))$, then
%
\begin{eqnarray}\label{f3int}
\lim_{n\to\infty}\mathop{\sum_{j=1}}_{j\ \mathrm{odd}}^{\lfloor
nt\rfloor}
g(F(t_{j-1}),t_{j-1})\Delta F_j^3
&=& \lim_{n\to\infty}\mathop{\sum_{j=1}}_{j\ \mathrm{even}}^{\lfloor
nt\rfloor}
g(F(t_{j-1}),t_{j-1})\Delta F_j^3\nonumber\\[-8pt]\\[-8pt]
&=& -\frac3{2\pi}\int_0^t g'(F(s),s) \,ds\nonumber
\end{eqnarray}
and
%
\begin{eqnarray}\label{b3int}
\lim_{n\to\infty}\mathop{\sum_{j=1}}_{j\ \mathrm{odd}}^{\lfloor
nt\rfloor}
g(F(t_j),t_j)\Delta F_j^3
&=& \lim_{n\to\infty}\mathop{\sum_{j=1}}_{j\ \mathrm{even}}^{\lfloor
nt\rfloor}
g(F(t_j),t_j)\Delta F_j^3\nonumber\\[-8pt]\\[-8pt]
&=& \frac3{2\pi}\int_0^t g'(F(s),s) \,ds,\nonumber
\end{eqnarray}
where the convergence is ucp.
\end{theorem}
\begin{remark}
The nonzero limits result from the dependence between
$F(t_{j-1})$ and $\Delta F_j$ in (\ref{f3int}), and $F(t_j)$ and
$\Delta
F_j$ in (\ref{b3int}). Note that
\[
E[F(t_{j-1})\Delta F_j] = \Delta t\partial_t\rho
(t_{j-1},t_{j-1}+\varepsilon)
\]
for some $0<\varepsilon<\Delta t$.
Similarly, $E[F(t_j)\Delta F_j]=\Delta t\partial_t\rho
(t_j,t_j-\varepsilon)$. If $X$
is a centered, quartic variation Gaussian process, then
\begin{eqnarray*}
\rho(s,t) &=& \tfrac12\bigl(EX(t)^2 + EX(s)^2 - E|X(t) - X(s)|^2\bigr)\\
&\approx& \tfrac12\bigl(EX(t)^2 + EX(s)^2 - |t - s|^{1/2}\bigr),
\end{eqnarray*}
which means the leading term in $\partial_t\rho(s,t)$ is
$-|t-s|^{-1/2}\operatorname{sgn}
(t-s)$. Hence, it is not surprising that the limits in (\ref{f3int})
and (\ref{b3int}) are of equal magnitude and opposite sign.
\end{remark}
\begin{pf*}{Proof of Theorem \ref{T:3ints}}
We prove only the case for odd indices. The proof for even
indices is nearly identical. To simplify notation, we will not
explicitly indicate that the indices are odd in the subscript of
the summation symbol (this convention applies only in this proof).

Using the truncation argument in the proof of Theorem
\ref{T:expan1}, we may assume that $g$ has compact support. Fix
$T>0$. Let $0\le s\le t\le T$ be arbitrary. Recall $\sigma_j$ and
$\widehat\sigma_j$ from Section \ref{sec:kb1}. Let
\begin{eqnarray*}
Z_n(t) &=& \sum_{j=1}^{\lfloor nt\rfloor} g(F(t_{j-1}),t_{j-1})\Delta
F_j^3,\\[-2pt]
X_n &=& X_n(s,t) = Z_n(t) - Z_n(s),\\
Y_n &=& Y_n(s,t) = 3\sum_{j=\lfloor ns\rfloor+1}^{\lfloor nt\rfloor}
g'(F(t_{j-1}),t_{j-1})\widehat\sigma_j\sigma_j^2.
\end{eqnarray*}
We may write
%
\begin{eqnarray}
&&
E|X_n - Y_n|^2\nonumber\\
&&\qquad= E \Biggl|\sum_{j=\lfloor ns\rfloor+1}^{\lfloor nt\rfloor}
g(F(t_{j-1}),t_{j-1})\Delta F_j^3\nonumber\\[-8pt]\\[-8pt]
&&\qquad\hspace*{20.7pt}{}
- 3\sum_{j=\lfloor ns\rfloor+1}^{\lfloor nt\rfloor}g'(F(t_{j-1}),t_{j-1})
\widehat\sigma_j\sigma_j^2
\Biggr|^2\nonumber\\
&&\qquad= (S_1 - S_2) - (S_2 - S_3),\nonumber
\end{eqnarray}
where
\begin{eqnarray*}
S_1 &=& \sum_{i=\lfloor ns\rfloor+1}^{\lfloor nt\rfloor}\sum
_{j=\lfloor ns\rfloor+1}^{\lfloor nt\rfloor}
E[g(F(t_{i-1}),t_{
i
-1})\Delta F_i^3 g(F(t_{j-1}),t_{j-1})\Delta F_j^3],\\
S_2 &=& 3\sum_{i=\lfloor ns\rfloor+1}^{\lfloor nt\rfloor}\sum
_{j=\lfloor ns\rfloor+1}^{\lfloor nt\rfloor}
E[g(F(t_{i-1}),t_{
i
-1})\Delta F_i^3 g'(F(t_{j-1}),t_{j-1})]\widehat\sigma_j\sigma
_j^2,\\
S_3 &=& 9\sum_{i=\lfloor ns\rfloor+1}^{\lfloor nt\rfloor}\sum
_{j=\lfloor ns\rfloor+1}^{\lfloor nt\rfloor}
E[g'(F(t_{i-1}),t_{
i
-1})g'(F(t_{j-1}),t_{j-1})]
\widehat\sigma_i\sigma_i^2\widehat\sigma_j\sigma_j^2.
\end{eqnarray*}
Let $\xi_1=F(t_{i-1})$, $\xi_2=\sigma_i^{-1}\Delta F_i$, $\xi
_3=F(t_{j-1})$,
$Y=\sigma_j^{-1}\Delta F_j$ and $\rho_k=E[\xi_kY]$. Define $f\in
C^3(\mathbb{R}^3)$
by $f(x)=g(x_1,t_{
i
-1})x_2^3g(x_3,t_{j-1})$ and define $h(x)=x^3$. By Corollary \ref
{C:GaussTaylor} with $k=2$,
\begin{eqnarray*}
E[f(\xi)Y^3] &=& \sum_{|\alpha|\le2}\frac1{\alpha!}\rho^\alpha
E[\partial^\alpha
f(\xi)]
E\bigl[h_{|\alpha|}(Y)Y^3\bigr] + R\\
&=& 3\sum_{|\alpha|=1}\frac1{\alpha!}\rho^\alpha E[\partial
^\alpha f(\xi)]
+ R,
\end{eqnarray*}
where $|R|\le C|\rho|^3$. Hence,
\begin{eqnarray*}
|S_1 - S_2|
&=& \Biggl|\sum_{i=\lfloor ns\rfloor+1}^{\lfloor nt\rfloor}\sum
_{j=\lfloor ns\rfloor+1}^{\lfloor nt\rfloor}
\sigma_i^3\sigma_j^3\bigl(E[f(\xi)Y^3] - 3\rho_3E[\partial_3f(\xi
)]\bigr) \Biggr|\\
&\le& C\sum_{i=\lfloor ns\rfloor+1}^{\lfloor nt\rfloor}\sum
_{j=\lfloor ns\rfloor+1}^{\lfloor nt\rfloor}
\sigma_i^3\sigma_j^3 (|\rho_1| + |\rho_2| + |\rho_3|^3)\\
&\le& C\sum_{i=\lfloor ns\rfloor+1}^{\lfloor nt\rfloor}\sum
_{j=\lfloor ns\rfloor+1}^{\lfloor nt\rfloor}
\bigl(\Delta t^{5/4}|EF(t_{i-1})\Delta F_j|\\
&&\hspace*{87pt}{} + \Delta t|E\Delta
F_i\Delta F_j|
+ \Delta t^{3/4}|\widehat\sigma_j|^3\bigr).
\end{eqnarray*}
By (\ref{cross}), (\ref{sighat}), Lemma \ref{L:sigdel}(i) with
$c=0$ and Lemma \ref{L:sigdel}(iii),
\begin{eqnarray*}
|S_1 - S_2|
&\le& C\sum_{i=\lfloor ns\rfloor+1}^{\lfloor nt\rfloor}\sum
_{j=\lfloor ns\rfloor+1}^{\lfloor nt\rfloor}
\bigl(\Delta t^{7/4}(|j - i|\vee1)^{-1/2}\\
&&\hspace*{87.4pt}{} + \Delta t^{3/2}(|j - i|\vee1)^{-3/2} + \Delta t^{9/4}\bigr)\\
&\le& C\sum_{i=\lfloor ns\rfloor+1}^{\lfloor nt\rfloor} \Delta t^{5/4}
\le C \biggl(\frac{\lfloor nt\rfloor- \lfloor ns\rfloor}n
\biggr)\Delta t^{1/4}.
\end{eqnarray*}
To estimate $S_2-S_3$, let $\xi_1=F(t_{i-1})$, $\xi_2=F(t_{j-1})$,
$Y=\sigma_i^{-1} \Delta F_i$ and $\rho_k=E[\xi_kY]$. Define $f\in
C^3(\mathbb{R}
^2)$ by
$f(x)=g(x_1,t_{j-1})g'(x_2,t_{j-1})$ and $h(x)=x^3$. As above,
\begin{eqnarray*}
|S_2 - S_3| &=& 3 \Biggl|
\sum_{i=\lfloor ns\rfloor+1}^{\lfloor nt\rfloor}\sum_{j=\lfloor
ns\rfloor+1}^{\lfloor nt\rfloor}
\widehat\sigma_j\sigma_j^2\sigma_i^3\bigl(E[f(\xi)Y^3] - 3\rho
_1E[\partial_1f(\xi
)]\bigr) \Biggr|\\
&\le& C\sum_{i=\lfloor ns\rfloor+1}^{\lfloor nt\rfloor}\sum
_{j=\lfloor ns\rfloor+1}^{\lfloor nt\rfloor}
|\widehat\sigma_j|\sigma_j^2\sigma_i^3(|\rho_2| + |\rho
_1|^3)\\
&\le& C\sum_{i=\lfloor ns\rfloor+1}^{\lfloor nt\rfloor}\sum
_{j=\lfloor ns\rfloor+1}^{\lfloor nt\rfloor}
\bigl(\Delta t^2(|j - i|\vee1)^{-1/2} + \Delta t^{5/2}\bigr)\\
&\le& C\sum_{i=\lfloor ns\rfloor+1}^{\lfloor nt\rfloor}\Delta t^{3/2}
\le C \biggl(\frac{\lfloor nt\rfloor- \lfloor ns\rfloor}n
\biggr)\Delta t^{1/2}.
\end{eqnarray*}
Combining these results, we have
\[
E|X_n - Y_n|^2
\le C \biggl(\frac{\lfloor nt\rfloor- \lfloor ns\rfloor}n
\biggr)\Delta t^{1/4}
\le C \biggl(\frac{\lfloor nt\rfloor- \lfloor ns\rfloor}n \biggr)^{5/4}.
\]
Note that
\begin{eqnarray*}
EY_n^2
&\le& C\sum_{i=\lfloor ns\rfloor+1}^{\lfloor nt\rfloor}\sum
_{j=\lfloor ns\rfloor+1}^{\lfloor nt\rfloor}
|\widehat\sigma_i\sigma_i^2\widehat\sigma_j\sigma_j^2|\\
&\le& C\sum_{i=\lfloor ns\rfloor+1}^{\lfloor nt\rfloor}\sum
_{j=\lfloor ns\rfloor+1}^{\lfloor nt\rfloor}
\Delta t^2
= C \biggl(\frac{\lfloor nt\rfloor- \lfloor ns\rfloor}n \biggr)^2.
\end{eqnarray*}
Since $t-s\le T$, this shows that
\[
E|Z_n(t) - Z_n(s)|^2 = EX_n^2
\le C \biggl(\frac{\lfloor nt\rfloor- \lfloor ns\rfloor}n \biggr)^{5/4}.
\]
Taking $s=0$ verifies condition (iii) of Theorem \ref{T:momcrit}.
Hence, by Corollary \ref{C:momcrit}, $\{Z_n\}$ is relatively
compact. Since $X_n - Y_n\to0$ in $L^2$, it will suffice, by Lemma
\ref{L:rcprob}, to show that
\[
Y_n(0,t) = 3\sum_{j=1}^{\lfloor nt\rfloor}
g'(F(t_{j-1}),t_{j-1})\widehat\sigma_j\sigma_j^2 \to-\frac
3{2\pi}\int_0^t
g'(F(s),s) \,ds
\]
in probability. For this, observe that by (\ref{sig2}) and (\ref{sighat}),
\begin{eqnarray*}
|\widehat\sigma_j\sigma_j^2 + \pi^{-1}\Delta t|
&\le&|\widehat\sigma_j + (2\pi)^{-1/2}\Delta t^{1/2}|\sigma
_j^2\\
&&{} + (2\pi)^{-1/2}\Delta t^{1/2}
|(2/\pi)^{1/2}\Delta t^{1/2} - \sigma_j^2|\\
&\le& Cj^{-1/2}\Delta t.
\end{eqnarray*}
Hence,
\[
\Biggl|\sum_{j=1}^{\lfloor nt\rfloor}g'(F(t_{j-1}),t_{j-1})\widehat
\sigma_j\sigma_j^2
+ \frac1\pi\sum_{j=1}^{\lfloor nt\rfloor
}g'(F(t_{j-1}),t_{j-1})\Delta t \Biggr|
\le C\Delta t^{1/2} \to0.
\]
Since
\[
\mathop{\sum_{j=1}}_{j\ \mathrm{odd}}^{\lfloor nt\rfloor
}g'(F(t_{j-1}),t_{j-1})\Delta t
\to\frac12\int_0^t g'(F(s),s) \,ds
\]
almost surely, this completes the proof of (\ref{f3int}).

For (\ref{b3int}), note that we may use (\ref{c1k3}) with $K=\mathbb{R}$
and $j=0$ to obtain
\begin{eqnarray*}
&&g(F(t_j),t_j) - g(F(t_{j-1}),t_{j-1})\\
&&\qquad= \int_{t_{j-1}}^{t_j}\partial_t g(F(t_j),s) \,ds
+ g(F(t_j),t_{j-1}) - g(F(t_{j-1}),t_{j-1})\\
&&\qquad= g'(F(t_{j-1}),t_{j-1})\Delta F_j + R,
\end{eqnarray*}
where $|R|\le C(\Delta t + \Delta F_j^2)$. Hence,
%
\begin{eqnarray}
\sum_{j=1}^{\lfloor nt\rfloor}g(F(t_j),t_j)\Delta F_j^3
&=& \sum_{j=1}^{\lfloor nt\rfloor}g(F(t_{j-1}),t_{j-1})\Delta
F_j^3\nonumber\\[-8pt]\\[-8pt]
&&{} + \sum_{j=1}^{\lfloor nt\rfloor}g'(F(t_{j-1}),t_{j-1})\Delta F_j^4 +
\widetilde R,\nonumber
\end{eqnarray}
where $|\widetilde R|\to0$ ucp. Applying (\ref{f3int}) and Theorem
\ref{T:4ints} completes the proof.
\end{pf*}

As a reminder, $X_n(t)\approx Y_n(t)$ means that $X_n-Y_n\to0$
ucp. Let
%
\begin{equation}\label{Jriem}
J_n(g,t) = \sum_{j=1}^{2\lfloor nt/2\rfloor}
g(F(t_{j-1}),t_{j-1})\Delta F_j^2(-1)^j.
\end{equation}
\begin{corollary}\label{C:expan2}
If $g\in C^{7,1}_3(\mathbb{R}\times[0,\infty))$, then
\[
I_n(g',t) \approx g(F(t),t) - g(F(0),0)
- \int_0^t \partial_tg(F(s),s) \,ds - \frac12J_n(g'',t),
\]
where $I_n(g,t)$ and $J_n(g,t)$ are given by (\ref{riem}) and
(\ref{Jriem}), respectively. Moreover,
\[
T_n^F(g',t) \approx g(F(t),t) - g(F(0),0)
- \int_0^t \partial_tg(F(s),s) \,ds,
\]
where $T_n^F$ is given by (\ref{Triem}).
\end{corollary}
\begin{pf}
By Theorems \ref{T:expan1}, \ref{T:trapexpan} and
\ref{T:3ints}, it will suffice to show that
\[
\sum_{j=1}^{\lfloor nt/2\rfloor}
g''(F(t_{2j-1}),t_{2j-1})(\Delta F_{2j}^2 - \Delta F_{2j-1}^2)
\approx J_n(g'',t).
\]
As before, we may assume that $g$ has compact support. Note that
\begin{eqnarray*}
&&\sum_{j=1}^{\lfloor nt/2\rfloor}
g''(F(t_{2j-1}),t_{2j-1})(\Delta F_{2j}^2 - \Delta F_{2j-1}^2)\\
&&\qquad= \mathop{\sum_{j=1}}_{j\ \mathrm{even}}^{2\lfloor nt/2\rfloor}
g''(F(t_{j-1}),t_{j-1})\Delta F_j^2
- \mathop{\sum_{j=1}}_{j\ \mathrm{odd}}^{2\lfloor nt/2\rfloor}
g''(F(t_j),t_j)\Delta F_j^2\\
&&\qquad= J_n(g'',t) - \mathop{\sum_{j=1}}_{j\ \mathrm{odd}}^{2\lfloor
nt/2\rfloor}
\{g''(F(t_j),t_j) - g''(F(t_{j-1}),t_{j-1})\}\Delta F_j^2.
\end{eqnarray*}
The proof is completed by using (\ref{kb5.2}) and applying Theorem
\ref
{T:3ints}.
\end{pf}
\begin{corollary}\label{C:expan3}
If $g\in C^{7,1}_3(\mathbb{R}\times[0,\infty))$, then
\[
T_n^X(g',t) \approx g(X(t),t) - g(X(0),0)
- \int_0^t \partial_tg(X(s),s) \,ds,
\]
where $T_n^X$ is given by (\ref{Triem}). This result remains true even
when $X=cF+\xi$, where $\xi$ satisfies only (\ref{exten0}), and is not
necessarily independent of $F$.
\end{corollary}
\begin{pf}
By passing to a subsequence, we may assume that Corollary \ref
{C:expan2} holds almost surely. It will therefore suffice to prove
Corollary \ref{C:expan3} under the assumption that $\xi$ is deterministic.

The claim is trivial when $c=0$. Suppose $c\ne0$. Let $h=h_\xi$ be
given by $h(x,t) =g(cx+\xi(t),t)$. We claim that $h\in
C^{7,1}_3(\mathbb{R}
\times[0,\infty))$. Note that $h^{(j)}(F(t),t)=c^jg^{(j)}(X(t),t)$ for
all $j\le7$. It is straightforward to verify (\ref{c1k0}) and (\ref
{c1k1}). Conditions (\ref{c1k2}) and (\ref{c1k3}) follow from the fact that
\[
\partial_t h^{(j)}(x,t) = c^jg^{(j+1)}\bigl(cx + \xi(t),t\bigr)\xi'(t)
+ c^j\partial_t g^{(j)}\bigl(cx + \xi(t),t\bigr)
\]
for all $j\le3$.

Observe that
\[
T_n^X(g',t) = T_n^F(h',t) + c^{-1}\sum_{j=1}^{\lfloor nt\rfloor}\frac{
h'(F(t_{j-1}),t_{j-1}) + h'(F(t_j),t_j)}2\Delta\xi_j.
\]
By our hypotheses on $\xi$, and the continuity of $h'$ and $F$, the
above summation converges to $\int_ 0^t h'(F(s),s)\xi'(s) \,ds$,
uniformly on compacts with probability one. Thus, by Corollary \ref
{C:expan2}, we have
\begin{eqnarray*}
T_n^X(g',t) &\approx& h(F(t),t) - h(F(0),0)\\
&&{}- \int_0^t \partial_th(F(s),s) \,ds
+ c^{-1}\int_0^t h'(F(s),s)\xi'(s) \,ds\\
&=& g(X(t),t) - g(X(0),0)
- \int_0^t \partial_tg(X(s),s) \,ds,
\end{eqnarray*}
which completes the proof.
\end{pf}

\section{Relative compactness}\label{S:rc}

The main result of this section is Theorem \ref{T:rcmain} below, from
which the relative compactness of $\{J_n(g,\cdot)\}$ will follow as a
corollary. [Recall that $\{J_n(g,\cdot)\}$ is defined in (\ref
{Jriem}).] Later in Section \ref{S:conv2Ito}, we will again need Theorem
\ref{T:rcmain}, when we show that $J_n$ converges weakly to an ordinary
It\^o integral.
\begin{theorem}\label{T:rcmain}
Let $g \in C^{7,1}_2(\mathbb{R}\times[0,\infty))$ have compact
support. Fix
$T>0$ and let $c$ and $d$ be integers such that $0\le t_c<t_d\le
T$. Then,
\[
E \Biggl|\sum_{j=c+1}^d \{g(F(t_{j-1}),t_{j-1})
- g(F(t_c),t_c) \}\Delta F_j^2(-1)^j \Biggr|^2
\le C|t_d - t_c|^{3/2},
\]
where $C$ depends only on $g$ and $T$.
\end{theorem}

Consider the simple case $c=0$ and $g(x,t)=x$. In that case, the above
expectation is
\[
E \Biggl|\sum_{j=1}^d F(t_{j-1})\Delta F_j^2(-1)^j \Biggr|^2
= \sum_{i=1}^d\sum_{j=1}^d
E[F(t_{i-1})\Delta F_i^2 F(t_{j-1})\Delta F_j^2](-1)^{i+j}.
\]
Using Corollary \ref{C:GaussTaylor}, we can remove the $\Delta F^2$
factors from inside the expectation. The leading term in the resulting
expansion would be roughly
\begin{eqnarray*}
&&\Delta t\sum_{i=1}^d\sum_{j=1}^d E[F(t_{i-1})F(t_{j-1})](-1)^{i+j}\\
&&\qquad= \Delta t\sum_{i,j\ \mathrm{even}}
E\bigl[\bigl(F(t_{i-1}) - F(t_{i-2})\bigr)\bigl(F(t_{j-1}) - F(t_{j-2})\bigr)\bigr].
\end{eqnarray*}
We could now use (\ref{cross}) to analyze these expectations and prove
the theorem in this simple case.

If we are to follow this strategy, then we will need an estimate
analogous to (\ref{cross}) which applies to functions of $F$. The
estimate in (\ref{cross}) was originally arrived at through direct
computations with the covariance function. Unfortunately, such direct
computations are not tractable for a general function of $F$. There is,
however, an alternative derivation of (\ref{cross}). Specifically, if
we observe that $|\partial_{st}\rho(s,t)| \le C|t-s|^{-3/2}$, where
$\partial
_{st}$ is the mixed second partial derivative, then we may conclude
that $|E[\Delta F_i\Delta F_j]| \le C\Delta t^2 |t_j - t_i|^{-3/2}$.
Based on
these heuristics, we begin with the following.
\begin{lemma}\label{L:derivest0}
Let $X$ be a centered Gaussian process with continuous covariance
function $\rho(s,t)$ and define $V(t)=\rho(t,t)$. Suppose that
$\rho$ is a $C^2$ function away from the set
$\{s=0\}\cup\{t=0\}\cup\{s=t\}$ and that $V(t)$ is a positive
$C^1$ function on $\{t>0\}$. Suppose that $\varphi\in C^2(\mathbb
{R})$ has
polynomial growth of order 2 with constants $K$ and $r$, and
define $V_\varphi(t) = E[\varphi(X(t))]$. Then,
\[
V_\varphi'(t) = \tfrac12 V'(t)E[\varphi''(X(t))].
\]
In particular, $|V_\varphi'(t)|\le C|V'(t)|$ for all $0<t\le T$, where
$C$ depends only on $K$, $r$ and $T$.
\end{lemma}
\begin{pf}
Let $\sigma(t)=V(t)^{1/2}$ and note that $\sigma$ is a positive
$C^1$ function on $\{t>0\}$. Fix $t>0$ and let
$X=\sigma(t)^{-1}X(t)$ so that $X$ is a standard normal random
variable and $V_\varphi(t)=E[\varphi(\sigma(t)X)]$. Since
$\varphi'$ has
polynomial growth, we may differentiate under the expectation,
giving
\[
V_\varphi'(t) = \sigma'(t) E[X \varphi'(\sigma(t)X)]
= \frac{V'(t)}{2\sigma(t)} E[\varphi'(\sigma(t)X)h_1(X)],
\]
where $h_n$ is given by (\ref{hermite}). By (\ref{hrmred2}), we have
\[
V'_\varphi(t) = \frac{V'(t)}{2\sigma(t)}E[\sigma(t)\varphi
''(\sigma(t)X)h_0(X)]
= \frac12 V'(t)E[\varphi''(X(t))].
\]
\upqed\end{pf}
\begin{prop}\label{P:derivest0}
Let $X$, $\rho$, and $V$ be as in Lemma \ref{L:derivest0}. Let
$g,h\in C^2(\mathbb{R})$ have polynomial growth of order 2 with common
constants $K$ and $r$, and define $f(s,t)= E[g(X(s)) h(X(t))]$.
Then,
%
\begin{eqnarray}
\label{derivest0.3}\hspace*{32pt}
\partial_sf(s,t) &=& \tfrac12
V'(s)E[g''(X(s))h(X(t))]\nonumber\\[-8pt]\\[-8pt]
&&{}
+ \partial_s\rho(s,t)E[g'(X(s))h'(X(t))]\quad\mbox{and}\nonumber\\
\label{derivest0.4}
\partial_tf(s,t) &=& \tfrac12 V'(t)E[g(X(s))h''(X(t))]
+ \partial_t\rho(s,t)E[g'(X(s))h'(X(t))]
\end{eqnarray}
whenever $0<s,t\le T$ and $s\ne t$. In particular,
\[
|\partial_s f(s,t)| \le C\bigl(|V'(s)| + |\partial_s\rho(s,t)|\bigr)
\]
and
\[
|\partial_t f(s,t)| \le C\bigl(|V'(t)| + |\partial_t\rho(s,t)|\bigr),
\]
where $C$ depends only on $K$, $r$ and $T$.
\end{prop}
\begin{pf}
By symmetry, we only need to prove (\ref{derivest0.3}). Let
$\sigma
(t)=V(t) ^{1/2}$ and note that $\sigma$ is a positive $C^1$ function on
$\{
t>0\}$. Let $r=r(s,t)=\sigma(s)^{-1}\sigma(t)^{-1}\times\rho(s,t)$ and define
$X_r=\sigma(s)^{-1}X(s)$ and $Y_r=\sigma(t)^{-1}X(t)$. Note that
$X_r$ and
$Y_r$ are jointly normal with mean zero, variance one and $E[X_rY_r]=r$.

Let $\varphi$ be as in Lemma \ref{L:cool?}. Then $f(s,t)=\varphi
(\sigma(s),\sigma
(t),r(s,t))$. Hence, by Lemmas \ref{L:cool} and \ref{L:cool?},
\begin{eqnarray*}
\partial_s f(s,t) &=& \sigma'(s)\sigma(s)E[g''(X(s))h(X(t))]\\
&&{} + \sigma'(s)\sigma(t)r(s,t)E[g'(X(s))h'(X(t))]\\
&&{} + \partial_s r(s,t)\sigma(s)\sigma(t)E[g'(X(s))h'(X(t))].
\end{eqnarray*}
Note that $\sigma'(s)=V'(s)/(2\sigma(s))$ and
\[
\partial_s r(s,t) = \frac{\partial_s\rho(s,t)}{\sigma(s)\sigma(t)}
- \frac{\rho(s,t)}{\sigma(s)^2\sigma(t)}\sigma'(s)
= \frac{\partial_s\rho(s,t)}{\sigma(s)\sigma(t)}
- \frac{V'(s)r(s,t)}{2\sigma(s)^2}.
\]
Thus,
\begin{eqnarray*}
\partial_s f(s,t) &=& \tfrac12 V'(s)E[g''(X(s))h(X(t))]\\
&&{} + \tfrac12 V'(s)\sigma(s)^{-1}\sigma
(t)r(s,t)E[g'(X(s))h'(X(t))]\\
&&{} + \partial_s\rho(s,t)E[g'(X(s))h'(X(t))]\\
&&{} - \tfrac12 V'(s)\sigma(s)^{-1}\sigma
(t)r(s,t)E[g'(X(s))h'(X(t))]\\
&=& \tfrac12 V'(s)E[g''(X(s))h(X(t))] + \partial_s\rho
(s,t)E[g'(X(s))h'(X(t))].\quad
\end{eqnarray*}
\upqed\end{pf}
\begin{theorem}\label{T:derivest0}
Let $X$, $\rho$ and $V$ be as in Lemma \ref{L:derivest0}. Let
$g,h\in C^3 (\mathbb{R})$ have polynomial growth of order 3 with common
constants $K$ and $r$, and define $f(s,t)= E[g(X(s)) h(X(t))]$.
Then
\[
|\partial_{st}f(s,t)| \le C|\partial_{st}\rho(s,t)|
+ C\bigl(|V'(s)| + |\partial_s\rho(s,t)|\bigr)\bigl(|V'(t)| +
|\partial_t\rho(s,t)|\bigr),
\]
whenever $0<s,t\le T$ and $s\ne t$, where $C$ depends only on $K$,
$r$ and $T$.
\end{theorem}
\begin{pf}
By (\ref{derivest0.3}),
\begin{eqnarray*}
\partial_{st}f(s,t) &=& \tfrac12 V'(s)\partial_t \{
E[g''(X(s))h(X(t))] \}
+ \partial_s\rho(s,t)\partial_t \{E[g'(X(s))h'(X(t))] \}\\
&&{}
+ \partial_{st}\rho(s,t)E[g'(X(s))h'(X(t))].
\end{eqnarray*}
Applying (\ref{derivest0.4}), we have
\begin{eqnarray*}
\partial_{st}f(s,t) &=& \tfrac14 V'(s)V'(t)E[g''(X(s))h''(X(t))]\\
&&{} + \tfrac12 V'(s)\partial_t\rho(s,t)E[g'''(X(s))h'(X(t))]\\
&&{} + \tfrac12 V'(t)\partial_s\rho(s,t)E[g'(X(s))h'''(X(t))]\\
&&{} + \partial_s\rho(s,t)\partial_t\rho(s,t)E[g''(X(s))h''(X(t))]\\
&&{} + \partial_{st}\rho(s,t)E[g'(X(s))h'(X(t))]
\end{eqnarray*}
and the theorem now follows.
\end{pf}

From Theorem \ref{T:derivest0}, we immediately obtain the following corollary.
\begin{corollary}\label{C:derivest0}
Let $X$, $\rho$ and $V$ be as in Lemma \ref{L:derivest0}. Let
$g,h\in C^3(\mathbb{R})$ have polynomial growth of order 3 with common
constants $K$ and $r$, and define $f(s,t)= E[g(X(s))h(X(t))]$. If
\begin{eqnarray*}
|V'(t)| &\le& Ct^{-1/2},\\
|\partial_s\rho(s,t)| + |\partial_t\rho(s,t)| &\le& C\bigl(s^{-1/2} +
(t - s)^{-1/2}\bigr)
\end{eqnarray*}
and
\[
|\partial_{st}\rho(s,t)| \le C\bigl(s^{-3/2} + (t - s)^{-3/2}\bigr)
\]
for all $0<s<t\le T$, where $C$ depends on only $T$, then
\[
|\partial_{st}f(s,t)| \le C\bigl(s^{-3/2} + (t - s)^{-3/2}\bigr)
\]
for a (possibly different) constant $C$ that depends only on $K$,
$r$ and $T$.
\end{corollary}

With this corollary in place, we can now begin proving Theorem
\ref{T:rcmain}.
\begin{lemma}\label{L:incrmom}
Suppose $g\in C^{1,1}_0(\mathbb{R}\times[0,\infty))$ has compact support.
If $p>0$, then
\[
E|g(F(t),t) - g(F(s),s)|^p \le C|t - s|^{p/4}
\]
for all $0\le s,t\le T$, where $C$ depends only on $g$, $p$ and $T$.
\end{lemma}
\begin{pf}
We write
\begin{eqnarray*}
&&g(F(t),t) - g(F(s),s)\\
&&\qquad= \int_s^t \partial_tg(F(t),u) \,du\\
&&\qquad\quad{}+ \bigl(F(t) - F(s)\bigr)\int_0^1 g'\bigl(F(s) + u\bigl(F(t) - F(s)\bigr),s\bigr) \,du.
\end{eqnarray*}
Hence, $|g(F(t),t) - g(F(s),s)| \le C|t - s| + C|F(t) - F(s)|$.
Since $F$ is a Gaussian process, an application of (\ref{sig1})
completes the proof.
\end{pf}
\begin{lemma}\label{L:claim01}
Recall that $\sigma_j^2=E\Delta F_j^2$. Under the hypotheses of Theorem
\ref{T:rcmain},
\[
E \Biggl|\sum_{j=c+1}^d \{g(F(t_{j-1}),t_{j-1})
- g(F(t_c),t_c) \}
\sigma_j^2(-1)^j \Biggr|^2 \le C|t_d - t_c|^{3/2},
\]
where $C$ depends only on $g$ and $T$.
\end{lemma}
\begin{pf}
By (\ref{sig2}),
\[
\sum_{j=c+1}^d \{g(F(t_{j-1}),t_{j-1})
- g(F(t_c),t_c) \}
\sigma_j^2(-1)^j = S + \varepsilon,
\]
where
\[
S = \biggl(\frac2\pi\biggr)^{1/2}\Delta t^{1/2}
\sum_{j=c+1}^d
\{g(F(t_{j-1}),t_{j-1})
- g(F(t_c),t_c) \}(-1)^j
\]
and, by H\"older's inequality,
\begin{eqnarray*}
|\varepsilon|^2 &\le& C\Delta t \Biggl(\sum_{j=c+1}^d
|g(F(t_{j-1}),t_{j-1}) - g(F(t_c),t_c)|j^{-3/2} \Biggr)^2\\
&\le& C\Delta t \Biggl(\sum_{j=c+1}^d
|g(F(t_{j-1}),t_{j-1}) - g(F(t_c),t_c)|^2 \Biggr) \Biggl(
\sum_{j=c+1}^d j^{-3} \Biggr).
\end{eqnarray*}
Hence, by Lemma \ref{L:incrmom},
\[
E|\varepsilon|^2 \le C\Delta t^{3/2}\sum_{j=c+1}^d |j-c|^{1/2}
\le C\Delta t^{3/2}|d-c|^{3/2} = C|t_d - t_c|^{3/2}.
\]
As for $S$, we assume, without loss of generality, that $c$ and
$d$ are both even. In that case,
\begin{eqnarray*}
S &=& \biggl(\frac2\pi\biggr)^{1/2}\Delta t^{1/2}
\mathop{\sum_{j=c+1}}_{j\ \mathrm{even}}^d
\{g(F(t_{j-1}),t_{j-1}) - g(F(t_{j-2}),t_{j-2}) \}\\[-1.5pt]
&=& \biggl(\frac2\pi\biggr)^{1/2}\Delta t^{1/2}
\mathop{\sum_{j=c+1}}_{j\ \mathrm{even}}^d
\biggl\{\int_{t_{j-2}}^{t_{j-1}}\partial_tg(F(t_{j-1}),s) \,ds\\[-1.5pt]
&&\hspace*{94.3pt}{} + g(F(t_{j-1}),t_{j-2}) - g(F(t_{j-2}),t_{j-2}) \biggr\}.
\end{eqnarray*}
Using (\ref{c1k3}) with $j=0$, the integral is bounded by $C\Delta t$ and
we have $E|S|^2\le C\Delta t(|t_d -t_c|^2 + S_1+S_2)$, where
\begin{eqnarray*}
S_1 &=& \mathop{\sum_{j=c+1}}_{j\ \mathrm{even}}^d
E|g(F(t_{j-1}),t_{j-2}) - g(F(t_{j-2},t_{j-2}))|^2,\\[-1.5pt]
S_2 &=& 2\mathop{\sum_{i=c+1}}_{i\ \mathrm{even}}^d
\mathop{\sum_{j=i+2}}_{j\ \mathrm{even}}^d |E [
\{g(F(t_{i-1}),t_{i-2}) - g(F(t_{i-2}),t_{i-2}) \}\\[-1.5pt]
&&\hspace*{71.6pt}{} \times
\{g(F(t_{j-1}),t_{j-2}) - g(F(t_{j-2}),t_{j-2}) \} ]
|\\[-1.5pt]
&=& 2\mathop{\sum_{i=c+1}}_{i\ \mathrm{even}}^d
\mathop{\sum_{j=i+2}}_{j\ \mathrm{even}}^d \biggl|
\int_{t_{i-2}}^{t_{i-1}}\int_{t_{j-2}}^{t_{j-1}}
\partial_{st}f_{ij}(s,t) \,dt \,ds \biggr|
\end{eqnarray*}
and $f_{ij}(s,t)=E[g(F(s),t_{i-2})g(F(t),t_{j-2})]$. Note that $F$
is a Gaussian process satisfying the conditions of Corollary
\ref{C:derivest0}. Hence,
\begin{eqnarray*}
S_2 &\le& C\mathop{\sum_{i=c+1}}_{i\ \mathrm{even}}^d
\mathop{\sum_{j=i+2}}_{j\ \mathrm{even}}^d
\int_{t_{i-2}}^{t_{i-1}}\int_{t_{j-2}}^{t_{j-1}}
\bigl(s^{-3/2} + (t - s)^{-3/2}\bigr) \,dt \,ds\\[-1.5pt]
&\le& C\Delta t^{1/2}\mathop{\sum_{i=c+1}}_{i\ \mathrm{even}}^d
\mathop{\sum_{j=i+2}}_{j\ \mathrm{even}}^d
\bigl((i-2)^{-3/2} + (j-i-1)^{-3/2}\bigr)\\[-1.5pt]
&\le& C\Delta t^{1/2}(d-c) = C\Delta t^{-1/2}|t_d - t_c|.
\end{eqnarray*}
By Lemma \ref{L:incrmom}, we also have $S_1\le C\Delta
t^{-1/2}|t_d-t_c|$. Hence,
\[
E|S|^2 \le C\Delta t^{1/2}|t_d - t_c|.
\]
Combined with
the estimate on $E|\varepsilon|^2$, this completes the proof.
\end{pf}
\begin{lemma}\label{L:claim02}
Let $\widehat\sigma_{c,j}=E[(F(t_{j-1})-F(t_c))\Delta F_j]$. Under the
hypotheses of Theorem \ref{T:rcmain},
\[
E \Biggl|\sum_{j=c+1}^d g''(F(t_{j-1}),t_{j-1})
\widehat\sigma_{c,j}^2(-1)^j \Biggr|^2 \le C\Delta t|t_d - t_c|,
\]
where $C$ depends only on $g$ and $T$.
\end{lemma}
\begin{pf}
By Lemma \ref{L:sigdel}(i) applied with $c=0$, and
(\ref{sighat}),
we have
\begin{eqnarray*}
|\widehat\sigma_{c,j} + (2\pi)^{-1/2} \Delta t^{1/2}|
&=& |\widehat\sigma_j - E[F(t_c)\Delta F_j] + (2\pi)^{-1/2}
\Delta t^{1/2}|\\
&\le& C\Delta t^{1/2}(j - c)^{-1/2}.
\end{eqnarray*}
Hence, by Lemma \ref{L:sigdel}(i),
\[
|\widehat\sigma_{c,j}^2 - (2\pi)^{-1}\Delta t|
\le C\Delta t^{1/2}|\widehat\sigma_{c,j} + (2\pi)^{-1/2} \Delta t^{1/2}|
\le C\Delta t(j - c)^{-1/2}.
\]
Therefore,
\[
\sum_{j=c+1}^d g''(F(t_{j-1}),t_{j-1})\widehat\sigma_{c,j}^2(-1)^j
= S + \varepsilon,
\]
where
\[
S = (2\pi)^{-1}\Delta t\sum_{j=c+1}^d g''(F(t_{j-1}),t_{j-1})(-1)^j
\]
and
\[
|\varepsilon|^2 \le C\Delta t^2 \Biggl(\sum_{j=c+1}^d
(j-c)^{-1/2} \Biggr)^2
\le C\Delta t^2(d - c) = C\Delta t|t_d - t_c|.
\]
The proof that $E|S|^2\le C\Delta t |t_d-t_c|$ is similar to that in
the proof of Lemma \ref{L:claim01}, except that we must use
(\ref{c1k3}) with $j=2$.
\end{pf}
\begin{lemma}\label{L:new_for_final}
Under the hypotheses of Theorem \ref{T:rcmain},
we have
\[
E \Biggl|\sum_{j=c+1}^d \{g(F(t_c),t_{j-1})
- g(F(t_c),t_c) \}
\Delta F_j^2(-1)^j \Biggr|^2 \le C|t_d - t_c|^3,
\]
where $C$ depends only on $g$ and $T$.
\end{lemma}
\begin{pf}
Let $Y(t)=g(F(t_c),t) - g(F(t_c),t_c)$ and note that
\begin{eqnarray*}
&&E \Biggl|\sum_{j=c+1}^d Y(t_{j-1})\Delta F_j^2(-1)^j \Biggr|^2\\
&&\qquad= \sum_{i=c+1}^d\sum_{j=c+1}^d
E[Y(t_{i-1})\Delta F_i^2 Y(t_{j-1})\Delta F_j^2](-1)^{i+j}.
\end{eqnarray*}
For fixed $i,j$, define $f\dvtx\mathbb{R}^2\to\mathbb{R}$ by
\[
f(x) = \biggl({\frac{g(x_1,t_{i-1}) - g(x_1,t_c)}{t_{i-1} - t_c}
} \biggr) \biggl({\frac{g(x_1,t_{j-1}) - g(x_1,t_c)}{t_{j-1} - t_c}
} \biggr)x_2^2.
\]
By (\ref{c1k3}) with $j=2$, $f$ has polynomial growth of order 2 with
constants $K$ and $r$ that do not depend on $i$ or $j$.

Let $\xi_1=F(t_c)$, $\xi_2=\sigma_i^{-1}\Delta F_i$, $Y=\sigma
_j^{-1}\Delta F_j$
and $h(y)=y^2$. By Corollary~\ref{C:GaussTaylor} with $k=1$, $E[f(\xi
)h(Y)] = E[f(\xi)]+R_1$, where $|R_1|\le C|\rho|^2$. Similarly, if
$\widetilde
f(x_1)=f(x_1,1)$, then
\[
E[f(\xi)] = E[\widetilde f(\xi_1)h(\xi_2)] = E[\widetilde f(\xi
_1)] + R_2,
\]
where $|R_2|\le C|E[\xi_1\xi_2]|^2$. Therefore,
\[
E[Y(t_{i-1})\Delta F_i^2 Y(t_{j-1})\Delta F_j^2]
= \sigma_i^2\sigma_j^2 E[Y(t_{i-1})Y(t_{j-1})] + R_3,
\]
where
\[
|R_3| = \sigma_i^2\sigma_j^2|t_{i-1} - t_c||t_{j-1} - t_c||R_1 + R_2|
\le\Delta t^3 |i - c||j - c||R_1 + R_2|.
\]
Using Lemma \ref{L:sigdel}(iii) and (\ref{cross}),
\begin{eqnarray*}
|\rho_1| &=& |E[\xi_1Y]| \le C\Delta t^{1/4}|j - c|^{-1/2},\\
|\rho_2| &=& |E[\xi_2Y]| \le C(|j - i|\vee1)^{-3/2},\qquad
|E[\xi_1\xi_2]| \le C\Delta t^{1/4}|i - c|^{-1/2}.
\end{eqnarray*}
This gives
\[
|R_3| \le C\Delta t^{7/2}(|i - c| + |j - c|)
+ C\Delta t^3(|j - i|\vee1)^{-3}(|i - c|^2 + |j - c|^2).
\]
Observe that
\[
\sum_{i=c+1}^d\sum_{j=c+1}^d |R_3(i,j)|
\le C\Delta t^{7/2}(d - c)^3 + C\Delta t^3(d - c)^3
\le C|t_d - t_c|^3.
\]
Hence, we are reduced to considering
\[
\sum_{i=c+1}^d\sum_{j=c+1}^d
\sigma_i^2\sigma_j^2E[Y(t_{i-1})Y(t_{j-1})](-1)^{i+j}
= E \Biggl|\sum_{j=c+1}^d Y(t_{j-1})\sigma_j^2(-1)^j \Biggr|^2.
\]
Using (\ref{c1k3}) with $j=0$ and (\ref{sig2}), we have
\begin{eqnarray*}
&&\Biggl|\sum_{j=c+1}^d Y(t_{j-1})\sigma_j^2(-1)^j \Biggr|\\
&&\qquad= \Biggl|\mathop{\sum_{j=c+1}}_{j\ \mathrm{even}}^d
\bigl(Y(t_{j-1})\sigma_j^2 - Y(t_{j-2})\sigma_{j-1}^2\bigr) \Biggr|\\
&&\qquad\le \mathop{\sum_{j=c+1}}_{j\ \mathrm{even}}^d
\bigl(|Y(t_{j-1})||\sigma_j^2 - \sigma_{j-1}^2|
 + |Y(t_{j-1}) - Y(t_{j-2})|\sigma_{j-1}^2\bigr)\\
&&\qquad\le C\sum_{j=c+1}^d (|t_{j-1} - t_c|j^{-3/2}\Delta t^{1/2}
+ \Delta t^{3/2})\\
&&\qquad\le C\Delta t^{3/2}\sum_{j=c+1}^d |j - c|^{-1/2}
\le C|t_d - t_c|^{3/2},
\end{eqnarray*}
which completes the proof.
\end{pf}
\begin{pf*}{Proof of Theorem \ref{T:rcmain}}
By Lemma \ref{L:new_for_final}, it will suffice to show that
\[
E \Biggl|\sum_{j=c+1}^d \{g(F(t_{j-1}),t_{j-1})
- g(F(t_c),t_{j-1}) \}\Delta F_j^2(-1)^j \Biggr|^2
\le C|t_d - t_c|^{3/2}.
\]
For brevity, let $\mathcal{X}(t)=g(F(t),t) - g(F(t_c),t)$ and write
\begin{eqnarray*}
&&E \Biggl|\sum_{j=c+1}^d \mathcal{X}(t_{j-1})\Delta F_j^2(-1)^j
\Biggr|^2\\
&&\qquad= \sum_{i=c+1}^d\sum_{j=c+1}^d
E[\mathcal{X}(t_{i-1})\Delta F_i^2 \mathcal{X}(t_{j-1})\Delta
F_j^2](-1)^{i+j}.
\end{eqnarray*}
Recall that $\sigma_j^2=E\Delta F_j^2$. Let $\delta
_c(t)=F(t)-F(t_c)$. Let
$\sigma_{c,j}^2 =E\delta_c(t_j)^2$. Let $\xi_1= F(t_c)$,
$\xi_2=\sigma_{c,i-1}^{-1}\delta_c(t_{i-1})$, $\xi_3=\sigma_{c,j-1}^{-1}
\delta_c(t_{j-1})$, $\xi_4=\sigma_i^{-1}\Delta F_i$ and
$\xi=(\xi_1,\ldots,\xi_4)$. For $x\in\mathbb{R}^4$, define
$f=f_{ij}$ by
\begin{eqnarray*}
f(x) &=& \biggl(
\frac{g(x_1 + \sigma_{c,i-1}x_2,t_{i-1})
- g(x_1,t_{i-1})}{\sigma_{c,i-1}} \biggr)\\
&&{}\times\biggl(\frac{g(x_1 + \sigma_{c,j-1}x_3,t_{j-1})
- g(x_1,t_{j-1})}{\sigma_{c,j-1}} \biggr)x_4^2.
\end{eqnarray*}
Let $Y=\sigma_j^{-1}\Delta F_j$ and $h(y)=y^2$.

Note that for $\theta\in(0,1]$ and $t_j\in[0,T]$,
$x\mapsto\theta^{-1}(g(x_1 + \theta x_2,t_j) - g(x_1,t_j))$ has
polynomial growth of order 6 with constants $K$ and $r$ that do
not depend on $\theta$ or $j$. Hence, $f$ has polynomial growth of
order 6 with constants $K$ and $r$. Thus, by Corollary
\ref{C:GaussTaylor} with $k=5$, if $\sigma=\sigma_{c,i-1}
\sigma_{c,j-1}\sigma_i^2\sigma_j^2$, then
\begin{eqnarray*}
&&E[\mathcal{X}(t_{i-1})\Delta F_i^2 \mathcal{X}(t_{j-1})\Delta
F_j^2]\\
&&\qquad= \sigma E[f(\xi)h(Y)]\\
&&\qquad= \sigma \biggl(\sum_{|\alpha|\le5}\frac1{\alpha!}\rho
^{\alpha}
E[\partial^\alpha f(\xi)]E\bigl[h_{|\alpha|}(Y)Y^2\bigr] + R_1 \biggr),
\end{eqnarray*}
where $\rho_j=E[\xi_jY]$ and $|R_1|\le C|\rho|^6$. If $p$ is a
positive integer, then by (\ref{hrm3}), (\ref{hrmx^n}) and
(\ref{hrmNu}) with $r=1$,
%
\begin{equation}\label{star11.9}
E\bigl[h_{|\alpha|} (Y)Y^p\bigr]=0 \qquad\mbox{if $p-|\alpha|$ is odd or
$|\alpha|>p$}.
\end{equation}
Hence, since $E[h_2(Y)Y^2]=E[Y^4-Y^2]=2$,
\begin{eqnarray*}
E[\mathcal{X}(t_{i-1})\Delta F_i^2 \mathcal{X}(t_{j-1})\Delta F_j^2]
= \sigma E[f(\xi)] + \sigma\rho_3^2 E[\partial_3^2 f(\xi)] +
\sigma R_2,
\end{eqnarray*}
where $R_2$ incorporates all terms of the form $\rho^{\alpha}
E[\partial^\alpha f(\xi)]$ with $|\alpha| =2$, except $\alpha=
(0,0,2,0)$. It
follows that
%
\begin{eqnarray}\label{star12}
&&
E[\mathcal{X}(t_{i-1})\Delta F_i^2 \mathcal{X}(t_{j-1})\Delta
F_j^2] \nonumber\\
&&\qquad= \sigma_j^2 E[\mathcal{X}(t_{i-1})\Delta F_i^2 \mathcal{X}(t_{j-1})]
\\
&&\qquad\quad{}+ \widehat\sigma_{c,j}^2 E[\mathcal{X}(t_{i-1})\Delta F_i^2
g''(F(t_{j-1}),t_{j-1})]
+ \sigma R_2,\nonumber
\end{eqnarray}
where $\widehat\sigma_{c,j}=E[\delta_c(t_{j-1})\Delta F_j]$ and
\[
|R_2| \le C(|\rho_1|^2 + |\rho_2|^2 + |\rho_4|
+ |\rho_1\rho_2| + |\rho_1\rho_3| + |\rho_2\rho_3| + |\rho_3|^6).
\]
The terms $|\rho_1\rho_4|$, $|\rho_2\rho_4|$, $|\rho_3\rho_4|$
and $|\rho_4|^2$ are not listed on the right-hand side of the
above estimate because $|\rho_1\rho_4| + |\rho_2\rho_4| +
|\rho_3\rho_4| + |\rho_4|^2 \leq C |\rho_4|$. Using (\ref{sig1})
and Lemma \ref{L:sigdel}, we have
\begin{eqnarray*}
|\sigma| &\le& C\Delta t^{3/2}|i-c|^{1/4}|j-c|^{1/4},\\
|\rho_2| &\le& C|i-c|^{-1/4}(|j-i|^{-1/2} + |j-c|^{-1/2}),\\
|\rho_1| &\le& C\Delta t^{1/4}|j-c|^{-1/2}\\
&\le& C|i-c|^{-1/4}(|j-i|^{-1/2} + |j-c|^{-1/2}),\\
|\rho_3| &\le& C|j-c|^{-1/4},\\
|\rho_4| &\le& C|j-i|^{-3/2}.
\end{eqnarray*}
Note that the above factors of $|j-i|$ are actually $(|j-i|\vee1)$,
although we have omitted this to simplify the notation. These estimates
now yield
\begin{eqnarray*}
|\sigma R_2| &\le& C\Delta t^{3/2}(|i-c|^{-1/4}|j-c|^{1/4}|j-i|^{-1}\\
&&\hspace*{37.4pt}{} + |i-c|^{-1/4}|j-c|^{-3/4}
+ |i-c|^{1/4}|j-c|^{-5/4}\\
&&\hspace*{37.4pt}{}+ |i-c|^{-1/4}|j-c|^{1/4}|j-i|^{-3/2}
+ |j-i|^{-1/2} + |j-c|^{-1/2}).
\end{eqnarray*}
Using $|j-c|\le|j-i|+|i-c|$ and $|i-c|\le|j-i|+|j-c|$, we can
show that
\[
|\sigma R_2| \le C\Delta t^{3/2}(|i-c|^{1/4}|j-c|^{-5/4}
+ |j-i|^{-1/2} + |j-c|^{-1/2})
\]
and, therefore,
that
\[
\sum_{i=c+1}^d\sum_{j=c+1}^d |\sigma R_2| \le C\Delta t^{3/2}
\sum_{i=c+1}^d(d - c)^{1/2}\le C\Delta t^{3/2}(d - c)^{3/2}
= C|t_d - t_c|^{3/2}.
\]
By (\ref{star12}), we are now reduced to considering the sums
%
\begin{eqnarray}\label{star13}\qquad
&&
\sum_{i=c+1}^d\sum_{j=c+1}^d \sigma_j^2
E[\mathcal{X}(t_{i-1})\Delta F_i^2
\mathcal{X}(t_{j-1})](-1)^{i+j}\nonumber\\[-8pt]\\[-8pt]
&&\qquad{} + \sum_{i=c+1}^d\sum_{j=c+1}^d\widehat\sigma_{c,j}^2
E[\mathcal{X}(t_{i-1})\Delta F_i^2
g''(F(t_{j-1}),t_{j-1})](-1)^{i+j},\nonumber
\end{eqnarray}
which will require two more applications of Corollary
\ref{C:GaussTaylor}. We will be brief in our presentation because
the following estimates can be obtained in a way very similar to
the one presented above.

For $x\in\mathbb{R}^3$, define $\widetilde f_1(x)=f(x_1,x_2,x_3,1)$.
Let $\widetilde
Y=\xi_4$, $\widetilde\xi=(\xi_1,\xi_2,\xi_3)$ and $\widetilde
\rho_j=E[\xi
_j\widetilde
Y]$. Note that $\widetilde f_1$ and $\widetilde f_2=\sigma
_{c,j-1}^{-1}\partial_3^2\widetilde
f$ both have polynomial growth of order 5 with constants $K$ and
$r$. Applying Corollary \ref{C:GaussTaylor} with $k=4$ and using
(\ref{star11.9}),
we have
%
\begin{eqnarray}\label{star14}
&&
\sigma_j^2 E[\mathcal{X}(t_{i-1})\Delta F_i^2 \mathcal
{X}(t_{j-1})]\nonumber\\
&&\qquad= \sigma E[\widetilde f_1(\widetilde\xi)h(\widetilde Y)]\nonumber\\
&&\qquad= \sigma E[\widetilde f_1(\widetilde\xi)]
+ \sigma\widetilde\rho_2^2 E[\partial_2^2\widetilde
f_1(\widetilde\xi)] + \sigma R_3\\
&&\qquad= \sigma_i^2\sigma_j^2 E[\mathcal{X}(t_{i-1})\mathcal
{X}(t_{j-1})]\nonumber\\
&&\qquad\quad{} + \widehat\sigma_{c,i}^2\sigma_j^2
E[g''(F(t_{i-1}),t_{i-1})\mathcal{X}(t_{j-1})]
+ \sigma R_3,\nonumber
\end{eqnarray}
where
%
\begin{equation}\label{star15}
|R_3| \le C(|\widetilde\rho_1|^2 + |\widetilde\rho_3|^2 +
|\widetilde\rho_1\widetilde\rho_2|
+ |\widetilde\rho_1\widetilde\rho_3| + |\widetilde\rho
_2\widetilde\rho_3| + |\widetilde\rho_2|^5).
\end{equation}
As before,
\begin{eqnarray*}
|\widetilde\rho_3| &\le& C|j-c|^{-1/4}(|i-j|^{-1/2} + |i-c|^{-1/2}),\\
|\widetilde\rho_1| &\le& C\Delta t^{1/4}|i-c|^{-1/2}\\
&\le& C|j-c|^{-1/4}(|i-j|^{-1/2} + |i-c|^{-1/2}),\\
|\widetilde\rho_2| &\le& C|i-c|^{-1/4},
\end{eqnarray*}
which gives
\[
|\sigma R_3| \le C\Delta t^{3/2}(|i-c|^{1/4}|i-j|^{-5/4}
+ |i-j|^{-1/2} + |i-c|^{-1/2})
\]
and shows that
%
\begin{equation}\label{star16}
\sum_{i=c+1}^d\sum_{j=c+1}^d |\sigma R_3|
\le C|t_d - t_c|^{3/2}.
\end{equation}
Similarly, if $\widetilde\sigma=\widehat\sigma_{c,j}^2\sigma
_i^2 \sigma_{c,i-1}$, then
%
\begin{eqnarray}\label{star17}
&&\widehat\sigma_{c,j}^2 E[\mathcal{X}(t_{i-1})\Delta F_i^2
g''(F(t_{j-1}), t_{j-1})]\nonumber\\
&&\qquad= \widetilde\sigma E[\widetilde f_2(\widetilde\xi)h(\widetilde
Y)]\nonumber\\
&&\qquad= \widetilde\sigma E[\widetilde f_2(\widetilde\xi)]
+ \widetilde\sigma\widetilde\rho_2^2E[\partial_2^2\widetilde
f_2(\widetilde\xi)] + \widetilde\sigma R_4\\
&&\qquad= \sigma_i^2\widehat\sigma_{c,j}^2 E[\mathcal
{X}(t_{i-1})g''(F(t_{j-1}),t_{j-1})]\nonumber\\
&&\qquad\quad{} + \widehat\sigma_{c,i}^2\widehat\sigma_{c,j}^2
E[g''(F(t_{i-1}),t_{i-1})g''(F(t_{j-1}),t_{i-1})]
+ \widetilde\sigma R_4,\nonumber
\end{eqnarray}
where $R_4$ also satisfies (\ref{star15}). Note that $|\widetilde
\sigma|\le
C\Delta
t^{7/4} |i-c|^{1/4}$. Since this is a better estimate than the one we
use for $|\sigma|$, the estimates above also give
%
\begin{equation}\label{star18}
\sum_{i=c+1}^d\sum_{j=c+1}^d |\widetilde\sigma R_4|
\le C|t_d - t_c|^{3/2}.
\end{equation}
By (\ref{star13}), (\ref{star14}), (\ref{star16}), (\ref{star17}) and
(\ref{star18}), we are reduced to considering the sums
\begin{eqnarray*}
&&
\sum_{i=c+1}^d\sum_{j=c+1}^d
\sigma_i^2\sigma_j^2 E[\mathcal{X}(t_{i-1})\mathcal
{X}(t_{j-1})](-1)^{i+j}\\
&&\qquad{} + \sum_{i=c+1}^d\sum_{j=c+1}^d
\widehat\sigma_{c,i}^2\sigma_j^2
E[g''(F(t_{i-1}),t_{i-1})\mathcal{X}
(t_{j-1})](-1)^{i+j}\\
&&\qquad{} + \sum_{i=c+1}^d\sum_{j=c+1}^d
\sigma_i^2\widehat\sigma_{c,j}^2 E[\mathcal{X}
(t_{i-1})g''(F(t_{j-1}),t_{j-1})](-1)^{i+j}\\
&&\qquad{} + \sum_{i=c+1}^d\sum_{j=c+1}^d\widehat\sigma_{c,i}^2\widehat
\sigma_{c,j}^2
E[g''(F(t_{i-1}),t_{i-1})g''(F(t_{j-1}),t_{j-1})](-1)^{i+j}.
\end{eqnarray*}
Note that this can be simplified to
\begin{eqnarray*}
&&
E \Biggl|\sum_{j=c+1}^d \sigma_j^2 \mathcal{X}(t_{j-1})(-1)^j
+ \sum_{j=c+1}^d
\widehat\sigma_{c,j}^2 g''(F(t_{j-1}),t_{j-1})(-1)^j \Biggr|^2\\
&&\qquad
\le C \Biggl(E \Biggl|
\sum_{j=c+1}^d \sigma_j^2 \mathcal{X}(t_{j-1})(-1)^j \Biggr|^2\\
&&\qquad\quad\hspace*{12.2pt}{} +
E \Biggl|
\sum_{j=c+1}^d
\widehat\sigma_{c,j}^2 g''(F(t_{j-1}),t_{j-1})(-1)^j \Biggr|^2 \Biggr).
\end{eqnarray*}
By Lemmas \ref{L:claim01} and \ref{L:claim02}, this completes the
proof.
\end{pf*}
\begin{corollary}\label{C:rcmain}
Recall $J_n(g,t)$ from (\ref{Jriem}). If $g\in
C^{7,1}_2(\mathbb{R}\times[0,\infty))$ has compact support, then
$\{J_n(g,\cdot)\}$ is relatively compact in $D_\mathbb{R}[0,\infty)$.
\end{corollary}
\begin{pf}
We shall apply Corollary \ref{C:momcrit} with $\beta= 4$. First,
note that $q(x+y)^4\le C(|x|^2+|y|^4)$. Fix $0\le s\le t\le T$.
Let $c=2\lfloor ns/2\rfloor$ and $d=2\lfloor nt/2\rfloor$. Then,
\begin{eqnarray*}
&&
E\bigl[q\bigl(J_n(t) - J_n(s)\bigr)^4\bigr]\\
&&\qquad\le CE \Biggl|\sum_{j=c+1}^d
\{g(F(t_{j-1}),t_{j-1}) - g(F(t_c),t_c) \}\Delta
F_j^2(-1)^j
\Biggr|^2\\
&&\qquad\quad{} + CE \Biggl|g(F(t_c),t_c)\sum_{j=c+1}^d\Delta F_j^2(-1)^j \Biggr|^4.
\end{eqnarray*}
By Theorem \ref{T:rcmain} and (\ref{Bmom}),
\[
E\bigl[q\bigl(J_n(t) - J_n(s)\bigr)^4\bigr] \le C|t_d-t_c|^{3/2} + C|t_d-t_c|^2
\le C \biggl(\frac{2\lfloor nt/2\rfloor-2\lfloor ns/2\rfloor}{n} \biggr)^{3/2}.
\]
This shows that one of the assumptions of Corollary
\ref{C:momcrit} holds. The other assumption follows from the same
estimate applied with $s=0$. By Corollary \ref{C:momcrit},
$\{J_n\}$ is relatively compact.
\end{pf}

\section{Convergence to a Brownian integral}\label{S:conv2Ito}

Recall that $J_n(g,t)$ is given by (\ref{Jriem}) and $B_n(t)$ is
given by (\ref{Bdef}). Note that
\[
J_n(g,t) = \kappa\int_0^t g(F_n(s-),N(s-)) \,dB_n(s),
\]
where $N(t)=\lfloor nt\rfloor/n$ and $F_n(t)=F( N(t))$.
In light of Theorem \ref{T:premain}, we would like to apply
Theorem \ref{T:KP}. Unfortunately, though, $\{B_n\}$ cannot be
decomposed in a way that satisfies (\ref{KPcond}). This is
essentially due to the numerous local oscillations of $B_n$. To
overcome this difficulty, we consider a modified version of $B_n$.

The process $B_n$ has a jump after every $\Delta t$ units of time. To
``smooth out'' this process, we shall restrict it so that it jumps only
after every $\Delta t^{1/4}$ units of time. Define
%
\begin{equation}\label{Bbar}
\overline B_n(t) = \kappa^{-1}\sum_{j=1}^{2m^3\lfloor mt/2\rfloor
}\Delta F_j^2(-1)^j,
\end{equation}
where $m=\lfloor n^{1/4}\rfloor$.
\begin{lemma}\label{L:Bbar}
The sequence $\{\overline B_n\}$ given by (\ref{Bbar}) satisfies
(\ref{KPcond}) and $B_n-\overline B_n\to0$ ucp.
\end{lemma}
\begin{pf}
Given $k$, let $d=d(k)=2m^3k$ and $c=c(k)=2m^3(k-1)$. Write
$\overline
B_n(t)= \kappa^{-1}\sum_{k=1}^{\lfloor mt/2\rfloor}\xi_k$, where
\[
\xi_k = \sum_{j=c+1}^d \Delta F_j^2(-1)^j.
\]
For $c<j\le d$, let $\Delta\overline F_j=\Delta F_j-E[\Delta
F_j|\mathcal{F}_{t_c}]$,
where $\mathcal{F}_t$ is given by (\ref{filt}). Let
\[
\overline\xi_k = \sum_{j=c+1}^d \Delta\overline F{}^2_j(-1)^j
\]
so that $\{\overline\xi_k\}$ is an i.i.d. sequence, by the remarks
following (\ref{filt}). In particular, $M_n(t)=\kappa^{-1}
\sum_{k=1}^{\lfloor mt/2\rfloor}\overline\xi_k$ is a martingale.
Let $A_n= \overline
B_n-M_n$. We must now verify (\ref{KPcond}).

Since $\{\Delta\overline F_j\}_{j=c+1}^\infty$ has the same law as
$\{\Delta F_j\}
_{j=1} ^\infty$, (\ref{Bmom}) implies
that
\[
E|\overline\xi_k|^2
= E \Biggl|\sum_{j=1}^{2m^3}\Delta F_j^2(-1)^j \Biggr|^2
= E|\kappa B_n(2m^3/n)|^2 \le Cn^{-1/4}.
\]
It follows that $E[M_n]_t = \kappa^{-1}\sum_{k=1}^{\lfloor
mt/2\rfloor}E|\overline\xi
_k|^2 \le Ct$ for all $n$. Also, by (\ref{sigcond2}),
\[
E|\xi_k - \overline\xi_k| \le C\Delta t^{1/2}\sum_{j=c+1}^d (j - c)^{-3/4}
\le C\Delta t^{1/2}(2m^3)^{1/4} \le Cn^{-5/16}.
\]
It follows that $EV_t(A_n)=\kappa^{-1}\sum_{k=1}^{\lfloor
mt/2\rfloor}E|\xi
_k-\overline\xi
_k| \le Ctn^{-1/16}$ and $\{\overline B_n\}$ satisfies (\ref{KPcond}).

By (\ref{Bmom}),
\[
E|\overline B_n(t) - \overline B_n(s)|^4 \le C \biggl(\frac
{2m^3\lfloor mt/2\rfloor- 2m^3\lfloor ms/2\rfloor}n \biggr)^2.
\]
By Corollary \ref{C:momcrit}, $\{\overline B_n\}$ is relatively
compact. By
Corollary \ref{C:conlim2} and Theorem \ref{T:premain}, $\{
B_n-\overline
B_n\}
$ is relatively compact. Hence, by Lemma \ref{L:rcprob}, in order to
show that $B_n-\overline B_n\to0$ ucp, it will suffice to show that
$B_n(t)-\overline B_n(t)\to0$ in probability for each fixed $t$.

For this, note that $n^{1/4}-1<m\le n^{1/4}$. Hence, $m^3\lfloor
mt/2\rfloor\le
nt/2$. Since $m^3\lfloor mt/2\rfloor$ is an integer, $m^3\lfloor
mt/2\rfloor\le\lfloor nt/2\rfloor$. By (\ref{Bmom}),
%
\begin{eqnarray}\label{star19}
E|B_n(t) - \overline B_n(t)|^4
&=& E \Biggl|\kappa^{-1}\sum_{j=2m^3\lfloor mt/2\rfloor+1}^{2\lfloor
nt/2\rfloor}
\Delta F_j^2(-1)^j \Biggr|^4\nonumber\\
&\le& C \biggl(\frac{2\lfloor nt/2\rfloor- 2m^3\lfloor mt/2\rfloor
}n \biggr)^2\nonumber\\[-8pt]\\[-8pt]
&\le& C \biggl(\frac{nt - m^4t + 2m^3}n \biggr)^2\nonumber\\
&\le& C \biggl(\frac{nt - (n^{1/4}-1)^4t + 2n^{3/4}}n \biggr)^2.\nonumber
\end{eqnarray}
Letting $n\to\infty$ completes the proof.
\end{pf}

With this lemma in place, we are finally ready to prove our main result.
\begin{theorem}\label{T:main}
Let $I_n(g,t)$ be given by (\ref{riem}) and $\kappa$, $B_n$ by
(\ref{kap}) and (\ref{Bdef}), respectively. Let $B$ be a standard
Brownian motion, independent of $F$. If $g\in C^{9,1}_4(\mathbb
{R}\times
[0,\infty))$, then $(F,B_n,I_n(g',\cdot)) \to(F,B,I^{F,B}(g',\cdot))$
in law in $D_{\mathbb{R}^3} [0,\infty)$, where $I^{F,B}(g',\cdot)$
is given by
(\ref{formint}).
\end{theorem}
\begin{remark}\label{R:append}
Suppose $\{W_n\}$ is another sequence of cadlag, $\mathbb{R}^
\ell
$-valued processes, adapted to a filtration of the form $\{\mathcal
{F}_t\vee
\mathcal{G}
_t^n\}$, where $\{\mathcal{F}_t\}$ and $\{\mathcal{G}_t^n\}$ are
independent. If
$(W_n,F,B_n)\to(W,F,B)$ in law in $D_{\mathbb{R}^ {\ell+2}}[0,\infty
)$, then
$(W_n,F,B_n,I_n(g',\cdot)) \to(W,F,B,I^{F,B}(g', \cdot))$ in law in
$D_{\mathbb{R}^{\ell+3}} [0,\infty)$. This can be seen by applying Remark
\ref
{R:KP} to (\ref{KPresult}) below.
\end{remark}
\begin{pf*}{Proof of Theorem \ref{T:main}}
By Lemma \ref{L:Bbar} and Theorem \ref{T:premain}, $\overline
B_n\to
B$ in law. Define $N(t)=2m^3\lfloor mt/2\rfloor/n$ and $\overline F_n(t)=F(
N(t))$. By continuity, $g''(\overline F_n(\cdot),N(\cdot))$
converges to
$g''(F(\cdot),\cdot)$ a.s. Hence, by Corollary \ref{C:conlim2} and
Lem\-ma~\ref{L:conlim3},
\[
(F, g''(\overline F_n (\cdot),N(\cdot)),\overline B_n)
\to(F,g''(F(\cdot),\cdot),B)
\]
in law in $D_{\mathbb{R}^3} [0,\infty)$.
Therefore, by Lemma \ref{L:Bbar}, Theorem \ref{T:KP} and Remark
\ref{R:KP},
%
\begin{eqnarray}\label{KPresult}
&&\biggl(F,g''(\overline F_n(\cdot),N(\cdot)),\overline B_n,
\kappa\int_0^\cdot g''(\overline F_n(s-),N(s-)) \,d\overline
B_n(s) \biggr)\nonumber\\[-8pt]\\[-8pt]
&&\qquad\to\biggl(F,g''(F(\cdot),\cdot),B,\kappa\int_0^\cdot g''(F(s),s)
\,dB(s) \biggr)\nonumber
\end{eqnarray}
in law in $D_{\mathbb{R}^4}[0,\infty)$. By Corollary \ref{C:expan2}
and Lemma
\ref{L:Bbar},
\begin{eqnarray*}
(F,B_n,I_n(g',t)) &\approx& \biggl(F, \overline B_n, g(F(\cdot),\cdot) -
g(F(0),0)- \int_0^t \partial_t g(F(s),s) \,ds\\
&&\hspace*{74.2pt}{}
- \frac\kappa2\int_0^\cdot g''(\overline F_n(s-),N(s-))
\,d\overline B_n(s) \biggr)\\
&&{} - \frac12 (0,0,\zeta_n(t)),
\end{eqnarray*}
where
\[
\zeta_n(t) = J_n(g'',t)
- \kappa\int_0^t g''(\overline F_n(s-),N(s-)) \,d\overline B_n(s).
\]
Hence, it will suffice to show that $\zeta_n\to0$ ucp.

By (\ref{Bbar}), $\overline B_n$ jumps only at times of the form
$s=2k/m$, where $k$ is an integer. At such a time,
$N(s-)=2m^3(k-1)/n$ and $\overline F_n(s-)=F( N(s-))$. Using the
notation in the proof of Lemma \ref{L:Bbar}, this gives
\begin{eqnarray*}
&&\kappa\int_0^t g''(\overline F_n(s-),N(s-)) \,d\overline B_n(s)\\
&&\qquad= \kappa\sum_{0<s\le t} g''(\overline F_n(s-),N(s-))\Delta
\overline B_n(s)\\
&&\qquad= \kappa\sum_{k=1}^{\lfloor mt/2\rfloor}
g''\bigl(F\bigl(t_{2m^3(k-1)}\bigr),t_{2m^3(k-1)}\bigr)
\kappa^{-1}\sum_{j=2m^3(k-1)+1}^{2m^3k}\Delta F_j^2(-1)^j\\
&&\qquad= \sum_{k=1}^{\lfloor mt/2\rfloor}\sum_{j=c+1}^d
g''(F(t_c),t_c)\Delta F_j^2(-1)^j.
\end{eqnarray*}
Hence, by (\ref{Jriem}), $\zeta_n(t) = \sum_{k=1}^{\lfloor
mt/2\rfloor}S_k +
\varepsilon
_n$, where
%
\begin{equation}\label{Sk}
S_k = \sum_{j=c+1}^d \{g''(F(t_{j-1}),t_{j-1})
- g''(F(t_c),t_c) \}\Delta F_j^2(-1)^j
\end{equation}
and
\[
\varepsilon_n = \sum_{j=2m^3\lfloor mt/2\rfloor+1}^{2\lfloor
nt/2\rfloor}
g''(F(t_{j-1}),t_{j-1})\Delta F_j^2(-1)^j.
\]
By the truncation argument in the proof of Theorem \ref{T:expan1},
we may assume that $g$ has compact support. Hence, by Corollary
\ref{C:rcmain}, $\{J_n(g'',\cdot)\}$ is relatively compact, so by
Corollary \ref{C:conlim2} and (\ref{KPresult}), $\{\zeta_n\}$ is
relatively compact. Therefore, by Lem\-ma~\ref{L:rcprob}, it will
suffice to show that $\zeta_n(t) \to0$ in probability for fixed $t$.

If $M=2m^3\lfloor mt/2\rfloor$ and $N=2\lfloor nt/2\rfloor$, then
\begin{eqnarray*}
\varepsilon_n &=& \sum_{j=M+1}^N \{g''(F(t_{j-1}),t_{j-1})
- g''(F(t_M),t_M) \}\Delta F_j^2(-1)^j\\
&&{} + g''(F(t_M),t_M)\sum_{j=M+1}^N\Delta F_j^2(-1)^j.
\end{eqnarray*}
Note that $g''$ is bounded and, by (\ref{Bdef}) and (\ref{Bmom}),
\[
E \Biggl|\sum_{j=M+1}^N\Delta F_j^2(-1)^j \Biggr|^4
= E|B_n(N/n) - B_n(M/n)|^4 \le C|t_N - t_M|^2.
\]
As in (\ref{star19}), this goes to zero as $n\to\infty$. Also, by
Theorem \ref{T:rcmain},
\[
E \Biggl|\sum_{j=M+1}^N \{g''(F(t_{j-1}),t_{j-1})
- g''(F(t_M),t_M) \}\Delta F_j^2(-1)^j \Biggr|^2
\le C|t_N-t_M|^{3/2}.
\]
Hence, $\varepsilon_n\to0$ in probability and it remains only to
check that
$\sum_ {k=1}^{\lfloor mt/2\rfloor} S_k\to0$ in probability.

Still using the notation from the proof of Lemma \ref{L:Bbar}, let
%
\begin{equation}\label{Skbar}
\overline S_k = \sum_{j=c+1}^d \{g''(F(t_{j-1}),t_{j-1})
- g''(F(t_c),t_c) \}\Delta\overline F{}^2_j(-1)^j,
\end{equation}
$\overline m_k=E[\overline S_k|\mathcal{F}_{t_c}]$ and $\overline
N_k=\overline S_k-\overline
m_k$. We claim that
%
\begin{equation}\label{marterr}
E|S_k - \overline N_k|^2 \le C\Delta t^{5/8}.
\end{equation}
For the moment, let us grant that this claim is true. In that
case,
\[
E \Biggl|\sum_{k=1}^{\lfloor mt/2\rfloor}S_k \Biggr|
\le\sum_{k=1}^{\lfloor mt/2\rfloor}E|S_k - \overline N_k|
+ \Biggl(E \Biggl|\sum_{k=1}^{\lfloor mt/2\rfloor}
\overline N_k \Biggr|^2 \Biggr)^{1/2}.
\]
Since $m\le n^{1/4}=\Delta t^{-1/4}$, (\ref{marterr}) gives
$\sum_{k=1}^{\lfloor mt/2\rfloor} E|S_k-\overline N_k|\le C\Delta
t^{1/16}\to0$. Also, if $k<\ell$, then $E[\overline N_k\overline
N_\ell] =
E[\overline N_k E[\overline N_\ell\mid\mathcal{F}_{t_{c(\ell)}}]]
= 0$. Hence,
\[
E \Biggl|\sum_{k=1}^{\lfloor mt/2\rfloor} \overline N_k \Biggr|^2
= \sum_{k=1}^{\lfloor mt/2\rfloor} E\overline N{}^2_k
\le C\sum_{k=1}^{\lfloor mt/2\rfloor}E|\overline N_k - S_k|^2
+ C\sum_{k=1}^{\lfloor mt/2\rfloor}ES_k^2.
\]
As above, the first summation goes to zero. For the second summation,
note that $g''\in C^{7,1}_2(\mathbb{R}\times[0,\infty))$ has compact support.
Thus, by (\ref{Sk}),
Theorem \ref{T:rcmain} and the fact that $d-c=2m^3\le2\Delta
t^{-3/4}$,
we have
\begin{eqnarray*}
ES_k^2 &=& E \Biggl|\sum_{j=c+1}^d \{g''(F(t_{j-1}),t_{j-1})
- g''(F(t_c),t_c) \}\Delta F_j^2(-1)^j \Biggr|^2\\
&\le& C|t_d - t_c|^{3/2} = C\Delta t^{3/2}(d-c)^{3/2} \le C\Delta t^{3/8}.
\end{eqnarray*}
Hence, $\sum_{k=1}^{\lfloor mt/2\rfloor}ES_k^2\le C\Delta t^{1/8}\to
0$, which
completes the proof of the theorem.

It remains only to prove (\ref{marterr}). By (\ref{Sk}) and
(\ref{Skbar}),
\begin{eqnarray*}
&& E|S_k - \overline S_k|^2\\
&&\qquad= E \Biggl|\sum_{j=c+1}^d
\{g''(F(t_{j-1}),t_{j-1}) - g''(F(t_c),t_c) \}
(\Delta F_j^2 - \Delta\overline F{}^2_j)(-1)^j \Biggr|^2\\
&&\qquad\le(d-c)\sum_{j=c+1}^d
E [{|g''(F(t_{j-1}),t_{j-1}) - g''(F(t_c),t_c)|^2
(\Delta F_j^2 - \Delta\overline F{}^2_j)^2} ].
\end{eqnarray*}
By H\"older's inequality, Lemma \ref{L:incrmom} and
(\ref{sigcond2}),
\begin{eqnarray*}
E|S_k - \overline S_k|^2 &\le& C(d - c)\sum_{j=c+1}^d
(t_j-t_c)^{1/2}\Delta t(j - c)^{-3/2}\\
&=& C\Delta t^{3/2}(d - c)\sum_{j=c+1}^d (j-c)^{-1}\\
&\le& C\Delta t^{3/2}(d - c)^{7/6}\le C\Delta t^{5/8}.
\end{eqnarray*}
Hence, it will suffice to show that $E|\overline m_k|^2 \le C\Delta
t^{5/8}$.

By (\ref{Skbar}),
\begin{eqnarray*}
\overline m_k &=& \sum_{j=c+1}^d
E[g''(F(t_{j-1}),t_{j-1})\Delta\overline F{}^2_j \mid\mathcal
{F}_{t_c}](-1)^j\\
&&{} - \sum_{j=c+1}^d g''(F(t_c),t_c)E[\Delta\overline F{}^2_j](-1)^j\\
&=& \sum_{j=c+1}^d
E\bigl[g''\bigl(G(t_{j-c-1}) + X_{j-c},t_{j-1}\bigr)\Delta G_{j-c}^2
\mid\mathcal{F}_{t_c}\bigr](-1)^j\\
&&{} - \sum_{j=c+1}^d g''(F(t_c),t_c)
E[\Delta G_{j-c}^2](-1)^j,
\end{eqnarray*}
where $G(t)=F(t+t_c)-E[F(t+t_c)\mid\mathcal{F}_{t_c}]$ and
$X_j=E[F(t_{j+c-1})\mid\mathcal{F}_ {t_c}]$. As noted in the discussion
following (\ref{filt}), $G$ is independent of $\mathcal{F}_{t_c}$ and has
the same law as $F$. Thus,
\[
\overline m_k = \mathop{\sum_{j=1}}_{j\ \mathrm{even}}^{d-c}
\bigl(\varphi_j(X_j) - \varphi_{j-1}(X_{j-1})\bigr)
- g''(F(t_c),t_c)\mathop{\sum_{j=1}}_{j\ \mathrm{even}}^{d-c}
(\sigma_j^2 - \sigma_{j-1}^2),
\]
where
\[
\varphi_j(x) = E\bigl[g''\bigl(F(t_{j-1}) + x,t_{j+c-1}\bigr)\Delta F_j^2\bigr].
\]
Using (\ref{hrmred2}), if $\sigma^2=EF(t_{j-1})^2$, then we have
\begin{eqnarray*}
\varphi_j(x) &=& \sigma_j^2E\bigl[g''\bigl(F(t_{j-1}) + x,t_{j+c-1}\bigr)\bigr]\\
&&{} + \sigma_j^2E\bigl[g''\bigl(\sigma(\sigma^{-1}F(t_{j-1})) + x,t_{j+c-1}\bigr)
h_2(\sigma_j^{-1}\Delta F_j)\bigr]\\
&=& \sigma_j^2E\bigl[g''\bigl(F(t_{j-1}) + x,t_{j+c-1}\bigr)\bigr]
+ \sigma_j^2(E[\sigma^{-1}F(t_{j-1})\sigma_j^{-1}\Delta F_j])^2\\
&&{} \times E\bigl[\sigma^2g^{(4)}\bigl(\sigma(\sigma^{-1}F(t_{j-1})) +
x,t_{j+c-1}\bigr)
h_0(\sigma_j^{-1}\Delta F_j)\bigr]\\
&=& \sigma_j^2E\bigl[g''\bigl(F(t_{j-1}) + x,t_{j+c-1}\bigr)\bigr]\\
&&{} + (E[F(t_{j-1})\Delta F_j])^2E\bigl[g^{(4)}\bigl(F(t_{j-1}) +
x,t_{j+c-1}\bigr)\bigr]\\
&=& \sigma_j^2 b_j(x) + \widehat\sigma_j^2 c_j(x),
\end{eqnarray*}
where
\begin{eqnarray*}
b_j(x) &=& E\bigl[g''\bigl(F(t_{j-1}) + x,t_{j+c-1}\bigr)\bigr],\\
c_j(x) &=& E\bigl[g^{(4)}\bigl(F(t_{j-1}) + x,t_{j+c-1}\bigr)\bigr].
\end{eqnarray*}
We may therefore write
%
\begin{equation}\label{mkbarest}
\overline m_k = \mathop{\sum_{j=1}}_{j\ \mathrm{even}}^{d-c}
\sum_{i=1}^5 \mathcal{E}_i,
\end{equation}
where
\begin{eqnarray*}
\mathcal{E}_1 &=& (\sigma_j^2 - \sigma^2_{j-1})b_j(X_j),\\
\mathcal{E}_2 &=& \sigma^2_{j-1}\bigl(b_j(X_j) - b_{j-1}(X_{j-1})\bigr),\\
\mathcal{E}_3 &=& (\widehat\sigma_j^2 - \widehat\sigma
_{j-1}^2)c_j(X_j),\\
\mathcal{E}_4 &=& \widehat\sigma_{j-1}^2\bigl(c_j(X_j) -
c_{j-1}(X_{j-1})\bigr),\\
\mathcal{E}_5 &=& - g''(F(t_c),t_c)(\sigma_j^2 - \sigma_{j-1}^2).
\end{eqnarray*}
For $\mathcal{E}_1$, (\ref{sig2}) gives $|\sigma^2_j - \sigma
^2_{j-1}| \le
Cj^{-3/2}\Delta t^{1/2}$. Hence,
\[
E \Biggl|\mathop{\sum_{j=1}}_{j\ \mathrm{even}}^{d-c}\mathcal
{E}_1 \Biggr|^2
\le C\Delta t.
\]
The same estimate also applies to $\mathcal{E}_5$. For $\mathcal
{E}_2$, let us write
\begin{eqnarray*}
&&|b_j(x_j) - b_{j-1}(x_{j-1})|\\
&&\qquad\le|b_j(x_j) - b_j(x_{j-1})|
+ |b_j(x_{j-1}) - b_{j-1}(x_{j-1})|\\
&&\qquad\le C|x_j - x_{j-1}|\\
&&\qquad\quad{} + \bigl|E\bigl[g''\bigl(F(t_{j-1}) + x_{j-1},t_{j+c-1}\bigr)
- g''\bigl(F(t_{j-2}) + x_{j-1},t_{j+c-2}\bigr)\bigr]\bigr|\\
&&\qquad\le C|x_j - x_{j-1}|\\
&&\qquad\quad{} + \bigl|E\bigl[g''\bigl(F(t_{j-1}) + x_{j-1},t_{j+c-1}\bigr)
- g''\bigl(F(t_{j-2}) + x_{j-1},t_{j+c-1}\bigr)\bigr]\bigr|\\
&&\qquad\quad{} + \bigl|E\bigl[g''\bigl(F(t_{j-2}) + x_{j-1},t_{j+c-1}\bigr)
- g''\bigl(F(t_{j-2}) + x_{j-1},t_{j+c-2}\bigr)\bigr]\bigr|\\
&&\qquad\le C|x_j - x_{j-1}| + |\beta_2'(t^*)|\Delta t
+ C\Delta t,
\end{eqnarray*}
where $\beta_2(t)=E[g''(F(t)+x_{j-1},t_{j+c-1})]$ and $t^*\in
(t_{j-2},t_{j-1})$, and where we have used (\ref{c1k3}) with $j=2$. By
Lemma \ref{L:derivest0}, $|\beta_2'(t)|\le Ct^{-1/2}$. Also, note that
$X_j-X_{j-1}=E[\Delta F_{j+c-1}\mid\mathcal{F}_{t_c}]$ so that by
(\ref{sigcond1}), $E|X_j-X_{j-1}|^2 \le Cj^{-3/2}\Delta t^{1/2}$. Thus,
%
\begin{eqnarray}\label{dagger1}
E \Biggl|\mathop{\sum_{j=1}}_{j\ \mathrm{even}}^{d-c}\mathcal
{E}_2 \Biggr|^2
&\le&\Biggl(\sum_{j=1}^{d-c}\sigma_{j-1}^4 \Biggr)
\Biggl(\sum_{j=1}^{d-c}E|b_j(X_j) - b_{j-1}(X_{j-1})|^2 \Biggr)\nonumber\\
&\le& C\Delta t^{1/4}\sum_{j=1}^{d-c}(j^{-3/2}\Delta t^{1/2}
+ j^{-1}\Delta t)\\
&\le& C\bigl(\Delta t^{3/4} + \Delta t^{5/4}(d - c)^{2/3}\bigr) \le C\Delta
t^{3/4}.\nonumber
\end{eqnarray}
For $\mathcal{E}_3$, (\ref{sighat}) gives $|\widehat\sigma
_j^2-\widehat\sigma_{j-1}^2|\le
Cj^{-1/2} \Delta t$. Hence,
\[
E \Biggl|\mathop{\sum_{j=1}}_{j\ \mathrm{even}}^{d-c}\mathcal
{E}_3 \Biggr|^2
\le C\Delta t^2(d-c) \le C\Delta t^{5/4}.
\]
For $\mathcal{E}_4$, as above, we have
\[
|c_j(x_j) - c_{j-1}(x_{j-1})|
\le C|x_j - x_{j-1}| + |\beta_4'(t^*)|\Delta t + C\Delta t,
\]
where $\beta_4(t)=E[g^{(4)}(F(t)+x_{j-1},t_{j+c-1})]$ and $t^*\in
(t_{j-2},t_{j-1})$, and where we have used (\ref{c1k3}) with $j=4$. It
therefore follows, as in (\ref{dagger1}), that
\[
E \Biggl|\mathop{\sum_{j=1}}_{j\ \mathrm{even}}^{d-c}\mathcal
{E}_4 \Biggr|^2
\le C\Delta t^{3/4}.
\]
Applying these five estimates to (\ref{mkbarest}) shows that
$E|\overline
m_k|^2 \le C\Delta t^{3/4} \le C\Delta t^{5/8}$ and completes the
proof.
\end{pf*}
\begin{corollary}\label{C:main}
Let $\xi$ be a
continuous
stochastic process, independent of $F$, such that (\ref{exten0})
holds. Let $X=cF+\xi$, where $c\in\mathbb{R}$. Let $I_n^X(g,t)$ be
given by
(\ref{riem}) and $\kappa$, $B_n$ by (\ref{kap}) and (\ref{Bdef}),
respectively. Let $B$ be a standard Brownian motion, independent of
$(F,\xi)$. If $g\in C^{9,1}_4(\mathbb{R}\times[0,\infty))$, then
$(F,\xi
,B_n,I_n^X(g',\cdot))\to(F,\xi,B,I^{X,c^2B} (g',\cdot))$ in law in
$D_{\mathbb{R}^4}[0,\infty)$, where $I^{X,Y}$ is given by (\ref{formint}).
\end{corollary}
\begin{remark}\label{R:main1}
Recall $Q_n^X$ from Section \ref{S:intro} and note that
$Q_n^F=\kappa B_n$. Note that $Q_n^X(t)\approx c^2Q_n^F(t)$ because
$\Delta X^2 = c^2 \Delta F^2 + o(\Delta t)$. This, together with Corollary
\ref{C:main}, implies that $(X,Q_n^X,I_n^X(g',\cdot))\to(X,\kappa
c^2B,I^{X,c^2B} (g',\cdot))$ in law in $D_{\mathbb{R}^3}[0,\infty)$.
\end{remark}
\begin{remark}\label{R:main2}
Suppose $\{W_n\}$ is another sequence of cadlag, $\mathbb{R}^
\ell
$-valued processes, adapted to a filtration of the form $\{\mathcal
{F}_t \vee
\mathcal{G}
_t^n\}$, where $\{\mathcal{F}_t\}$ and $\{\mathcal{G}_t^n\}$ are
independent. As in
Remark \ref{R:append}, if $(W_n,F,B_n)\to(W,F,B)$ in law in
$D_{\mathbb{R}^
{\ell+2}}[0,\infty)$, then $(W_n,F,\xi,B_n,I_n^X(g',\cdot)) \to
(W,F,\xi
,B,I^{X,c^2B}(g', \cdot))$ in law in $D_{\mathbb{R}^{\ell+4}}
[0,\infty)$.
\end{remark}
\begin{pf*}{Proof of Corollary \ref{C:main}}
The claim is trivial when $c=0$. Suppose $c\ne0$. We first assume
$\xi$ is deterministic. Let $h=h_\xi$ be given by $h(x,t) =g(cx+\xi
(t),t)$. We claim that $h\in C^{9,1}_4(\mathbb{R}\times[0,\infty))$.
Note that
$h^{(j)}(F(t),t)=c^jg^{(j)}(X(t),t)$ for all $j\le9$. It is
straightforward to verify (\ref{c1k0}) and (\ref{c1k1}). Conditions
(\ref{c1k2}) and (\ref{c1k3}) follow from the fact that
%
\begin{equation}\label{exten1}
\partial_t h^{(j)}(x,t) = c^jg^{(j+1)}\bigl(cx + \xi(t),t\bigr)\xi'(t)
+ c^j\partial_t g^{(j)}\bigl(cx + \xi(t),t\bigr)
\end{equation}
for all $j\le4$.

Observe that
\[
I_n^X(g',t) = I_n(h',t) + c^{-1}\sum_{j=1}^{\lfloor nt/2\rfloor}
h'(F(t_{2j-1}),t_{2j-1})\bigl(\xi(t_{2j}) - \xi(t_{2j-2})\bigr).
\]
By our hypotheses on $\xi$, and the continuity of $h'$ and $F$, the
above sum converges uniformly on compacts, with probability one, to
$\int_ 0^t h'(F(s),s)\xi'(s) \,ds$. Thus, by Theorem \ref{T:main} and
Remark \ref{R:append}, $(F,\xi,B_n,I_n^X (g',\cdot))\to(F,\xi
,B,\mathcal{I})$, where
\begin{eqnarray*}
\mathcal{I}&=& h(F(t),t) - h(F(0),0) - \int_0^t\partial_t h(F(s),s) \,ds
- \frac\kappa2\int_0^t h''(F(s),s) \,dB(s)\\
&&{} + c^{-1}\int_0^t h'(F(s),s)\xi'(s) \,ds.
\end{eqnarray*}
Using (\ref{exten1}) with $j=0$, this gives
\[
\mathcal{I}= g(X(t),t) - g(X(0),0) - \int_0^t\partial_t g(X(s),s) \,ds
- \frac{\kappa c^2}2\int_0^t g''(X(s),s) \,dB(s),
\]
completing the proof.

Now, suppose $\xi$ is random and independent of $F$. Let $H\dvtx D_{\mathbb{R}
^4}[0, \infty)\to\mathbb{R}$ be bounded and continuous. Since we
have proven
the result for deterministic $\xi$, it follows that
\[
E[H(F,\xi,B_n,I_n^X(g',\cdot)) \mid\xi]
\to E[H(F,\xi,B,I^{X,c^2B}(g',\cdot))
\mid\xi] \qquad\mbox{a.s.}
\]
Applying the dominated convergence theorem completes the proof.
\end{pf*}

We now give two examples of processes $X$ satisfying the
conditions of Corollary~\ref{C:main}.
\begin{example}\label{ex:1}
Consider the stochastic heat equation $\partial_t u=\frac
12\partial_x^2
u+\dot{W}(x,t)$ with initial conditions $u(x,0)=f(x)$. Under
suitable conditions on $f$, the unique solution is
\[
u(x,t) = \int_{\mathbb{R}\times[0,t]} p(x - y,t - r)W(dy\times dr)
+ v(t,x),
\]
where
\[
v(x,t) = \int_\mathbb{R}p(x - y,t)f(y) \,dy.
\]
For example, if $f$ has polynomial growth, then this is the unique
solution and, moreover, $\partial_tv$ is continuous on
$\mathbb{R}\times[0,\infty)$. This implies that $t\mapsto v(x,t)$
satisfies (\ref{exten0}). Hence, $X(t)=u(x,t)=F(t) +v(x,t)$
satisfies the conditions of Corollary~\ref{C:main}. This remains
true when $f$ is allowed to be a stochastic process, independent
of $W$.
\end{example}
\begin{example}\label{ex:2}
This example is based on a decomposition of bifractional Brownian
motion due to Lei and Nualart \cite{N2}. Let $W$ be a standard Brownian
motion, independent of $F$. Define
\[
\xi(t) = (16\pi)^{-1/4}\int_0^\infty(1 - e^{-st})s^{-3/4} \,dW(s).
\]
By Proposition 1 and Theorem 1 in \cite{N2},
we have $\xi\in C^1((0,\infty))$ a.s.
Moreover, if $c=(\pi/2)^{1/4}$, then $X=cF+\xi$ has the same law as
$B^{1/4}$, fractional Brownian motion with Hurst parameter $H=1/4$.
If $\varphi\in C^\infty[0,\infty)$ with $\varphi=0$ on
$[0,\varepsilon/4]$ and
$\varphi=1$
on $[\varepsilon/2,\infty)$, then $\varphi\xi$ satisfies (\ref
{exten0}) and we may
apply Corollary \ref{C:main} to $X_\varepsilon=cF+\varphi\xi$ to
obtain that
\begin{eqnarray*}
&&(X(t), Q_n^X(t), I_n^X(g',t))
- (X(\varepsilon), Q_n^X(\varepsilon), I_n^X(g',\varepsilon))\\
&&\qquad= \Biggl(X(t) - X(\varepsilon), \sum_{j=\lfloor n\varepsilon
/2\rfloor+1}^{\lfloor nt/2\rfloor}
(\Delta X_{2j}^2 - \Delta X_{2j-1}^2),\\
&&\hspace*{41.1pt}\sum_{j=\lfloor n\varepsilon/2\rfloor+1}^{\lfloor nt/2\rfloor}
g(X(t_{2j-1}),t_{2j-1})\bigl(X(t_{2j}) - X(t_{2j-2})\bigr) \Biggr)
\end{eqnarray*}
converges in law in $D_{\mathbb{R}^3}[\varepsilon,\infty)$ as $n\to
\infty$ to
\begin{eqnarray*}
&&(X(t), \kappa c^2B(t), I^{X,c^2B}(g',t))
- (X(\varepsilon), \kappa c^2B(\varepsilon),
I^{X,c^2B}(g',\varepsilon))\\
&&\qquad= \biggl(X(t) - X(\varepsilon), \kappa c^2\bigl(B(t)-B(\varepsilon)\bigr),
g(X(t),t) - g(X(\varepsilon),\varepsilon)\\
&&\qquad\quad\hspace*{27.5pt}{}- \int_\varepsilon^t \partial_t g(X(s),s) \,ds
- \frac{\kappa c^2}2\int_\varepsilon^t \partial_x^2 g(X(s),s)
\,dB(s) \biggr).
\end{eqnarray*}
\end{example}

%

%
\printaddresses

\end{document}